\documentclass[aos,MSNbibl,nameyear,dvips]{arximspdf}

\doi{10.1214/12-AOS1052} 
\volume{40}
\issue{6}
\pubyear{2012}
\firstpage{2823}
\lastpage{2849}

\makeatletter

\newtheorem{proposition}{Proposition}
\newtheorem{corollary}{Corollary}
\newtheorem{lemma}{Lemma}
\newtheorem{theorem}{Theorem}

\newproclaim{remark}{Remark}

\newcommand{\rright}{\right}
\newcommand{\lleft}{\left}
\newcommand{\rrVert}{\Vert}
\newcommand{\llVert}{\Vert}

\newcommand{\cal}{\mathcal}

\newcommand{\X}{{\mathsf{X}}}
\newcommand{\tr}{\operatorname{tr}}
\newcommand{\rank}{\operatorname{rank}}

\makeatother

\begin{document}
\begin{frontmatter}

\title{Convergence analysis of the Gibbs sampler for Bayesian general
linear mixed models with improper priors}
\runtitle{Convergence of the Gibbs sampler}

\begin{aug}
\author[A]{\fnms{Jorge Carlos} \snm{Rom\'{a}n}\corref{}\ead[label=e1]{jc.roman@vanderbilt.edu}}
\and
\author[B]{\fnms{James P.} \snm{Hobert}\thanksref{t1}\ead[label=e2]{jhobert@stat.ufl.edu}}
\runauthor{J. C. Rom\'{a}n and J. P. Hobert}
\affiliation{Vanderbilt University and University of Florida}
\address[A]{Department of Mathematics \\
Vanderbilt University \\
Nashville, Tennessee 37240 \\
USA \\
\printead{e1}} 
\address[B]{Department of Statistics \\
University of Florida \\
Gainesville, Florida 32611\\
USA \\
\printead{e2}}
\end{aug}

\thankstext{t1}{Supported by NSF Grants DMS-08-05860 and DMS-11-06395.}

\received{\smonth{4} \syear{2012}}
\revised{\smonth{9} \syear{2012}}

%
\begin{abstract}
Bayesian analysis of data from the general linear mixed model is
challenging because any nontrivial prior leads to an intractable
posterior density. However, if a conditionally conjugate prior
density is adopted, then there is a simple Gibbs sampler that can be
employed to explore the posterior density. A popular default among
the conditionally conjugate priors is an improper prior that takes a
product form with a flat prior on the regression parameter, and
so-called power priors on each of the variance components. In this
paper, a convergence rate analysis of the corresponding Gibbs
sampler is undertaken. The main result is a simple, easily-checked
sufficient condition for geometric ergodicity of the Gibbs--Markov
chain. This result is close to the best possible result in the
sense that the sufficient condition is only slightly stronger than
what is required to ensure posterior propriety. The theory
developed in this paper is extremely important from a practical
standpoint because it guarantees the existence of central limit
theorems that allow for the computation of valid asymptotic standard
errors for the estimates computed using the Gibbs sampler.
\end{abstract}

%

\end{frontmatter}

\section{Introduction}
\label{secintro}

The general linear mixed model (GLMM) takes the form
%
\begin{equation}
\label{eqGLMM} Y = X \beta+ Zu + e,
\end{equation}
where $Y$ is an $N \times1$ data vector, $X$ and $Z$ are known matrices
with dimensions $N \times p$ and $N \times q$, respectively, $\beta$ is
an unknown $p \times1$ vector of regression coefficients, $u$ is a
random vector whose elements represent the various levels of the random
factors in the model and $e \sim\mathrm{N}_N(0, \sigma^2_e I)$. The
random vectors $e$\vadjust{\goodbreak} and $u$ are assumed to be independent. Suppose there
are $r$ random factors in the model. Then $u$ and $Z$ are partitioned
accordingly as $u =  ( u_1^T\enskip u_2^T \enskip\cdots\enskip u_r^T
)^T$ and $Z = (Z_1\enskip Z_2\enskip \cdots\enskip Z_r)$, where $u_i$
is $q_i \times 1$, $Z_i$ is $N \times q_i$ and $q_1 + \cdots+ q_r = q$.
Then
\[
Zu = \sum_{i=1}^r Z_i
u_i,
\]
and it is assumed that $u \sim\mathrm{N}_q(0, D)$, where $D =
\bigoplus_{i=1}^r \sigma^2_{u_i} I_{q_i}$. Let $\sigma^2$ denote the
vector of variance components, that is, $\sigma^2 = (\sigma^2_e\enskip
\sigma^2_{u_1} \enskip\cdots\enskip \sigma^2_{u_r})^T$.
For\vspace*{1pt} background on this model, which is sometimes called
the \textit{variance components model}, see
\citet{searcasemccu1992}.

A Bayesian version of the GLMM can be assembled by specifying a prior
distribution for the unknown parameters, $\beta$ and $\sigma^2$. A
popular choice is the proper (conditionally) conjugate prior that
takes $\beta$ to be multivariate normal, and takes each of the
variance components to be inverted gamma. One obvious reason for
using such a prior is that the resulting posterior has conditional
densities with standard forms, and this facilitates the use of the
Gibbs sampler.

In situations where there is little prior information, the
hyperparameters of this proper prior are often set to extreme values
as this is thought to yield a ``noninformative'' prior.
Unfortunately, these extreme proper priors approximate improper priors
that correspond to improper posteriors, and this results in various
forms of instability. This problem has led several authors, including
\citet{Dani1999} and \citet{Gelm2006}, to discourage the use
of such
extreme proper priors, and to recommend alternative default priors
that are improper, but lead to proper posteriors. Consider, for
example, the one-way random effects model given by
%
\begin{equation}
\label{eqowm} Y_{ij} = \beta+ \alpha_i + e_{ij},
\end{equation}
where $i=1,\ldots,c$, $j=1,\ldots,n_i$, the $\alpha_i$'s are i.i.d.
$\mathrm{N}(0,\sigma^2_\alpha)$, and the $e_{ij}$'s, which are
independent of the $\alpha_i$'s, are i.i.d. $\mathrm{N}(0,\sigma_e^2)$.
This is an important special case of model (\ref{eqGLMM}). (See
Section~\ref{secowm} for a detailed explanation of how the GLMM
reduces to the one-way model.) The standard diffuse prior for this
model, which is among those recommended by \citet{Gelm2006}, has
density $1 /  ( \sigma^2_e \sqrt{\sigma^2_\alpha}  )$. This
prior, like many of the improper priors for the GLMM that have been
suggested and studied in the literature, is called a ``power prior''
because it is a product of terms, each a variance component brought to
a (possibly negative) power. Of course, like the proper conjugate
priors mentioned above, power priors also lead to posteriors whose
conditional densities have standard forms.

In this paper, we consider the following parametric family of priors
for $(\beta,\sigma^2)$:
%
\begin{equation}
\label{eqbeta}\quad
p\bigl(\beta,\sigma^2;a,b\bigr) = \bigl(
\sigma^2_e\bigr)^{-(a_e+1)} e^{-{b_e}/{\sigma^2_e}} \Biggl[
\prod_{i=1}^r \bigl(\sigma^2_{u_i}
\bigr)^{-(a_i+1)} e^{-{b_i}/{\sigma^2_{u_i}}} \Biggr] I_{\mathbb{R}^{r+1}_+}\bigl(
\sigma^2\bigr),
\end{equation}
where $a=(a_e,a_1,\ldots,a_r)$ and $b=(b_e,b_1,\ldots,b_r)$ are fixed
hyperparameters, and $\mathbb{R}_+:= (0,\infty)$. By taking $b$ to
be the vector of 0's, we can recover the power priors described above.
Note that $\beta$ does not appear on the right-hand side of
(\ref{eqbeta}); that is, we are using a so-called flat prior for
$\beta$. Consequently, even if all the elements of $a$ and $b$ are
strictly positive, so that every variance component gets a proper
prior, the overall prior remains improper. There have been several
studies concerning posterior propriety in this context, but it is
still not known exactly which values of $a$ and $b$ yield proper
posteriors. The best known result is due to \citet{SunTsutHe2001},
and we state it below so that it can be used in a comparison later in
this section.

Define $\theta= (\beta^T\enskip u^T)^T$ and $W = (X\enskip Z)$, so that $W
\theta= X \beta+ Zu$. Let $y$ denote the observed value of $Y$, and
let $\phi_d(x; \mu, \Sigma)$ denote the $\mathrm{N}_d(\mu,\Sigma)$
density evaluated at the vector $x$. By definition, the posterior
density is proper if
\[
m(y):= \int_{\mathbb{R}_+^{r+1}} \int_{\mathbb{R}^{p+q}} \pi^*\bigl(
\theta,\sigma^2|y\bigr) \,d\theta \,d\sigma^2 < \infty,
\]
where
%
\begin{equation}
\label{eqpistar} \pi^*\bigl(\theta,\sigma^2|y\bigr) =
\phi_N\bigl(y;W\theta,\sigma^2_{e}I\bigr)
\phi_q (u;0,D ) p\bigl(\beta, \sigma^2;a,b\bigr).
\end{equation}
A routine calculation shows that the posterior is improper if
$\rank(X) < p$. The following result provides sufficient (and nearly
necessary) conditions for propriety. (Throughout the paper, the
symbol $P$ subscripted with a matrix will denote the projection onto
the column space of that matrix.)
%
\begin{theorem}[{[\citet{SunTsutHe2001}]}]
\label{thmpropriety}
Assume that $\operatorname{rank}(X)=p$, and let $t = \operatorname
{rank} ( Z^T
(I-P_X) Z )$ and $\mathrm{SSE} = \llVert {(I-P_W) y}\rrVert^2$. If the
following four conditions hold, then $m(y) < \infty$:
\begin{longlist}[(D)]
\item[(A)] For each $i \in\{1,2,\ldots,r\}$, one of the following
holds:
\[
\mbox{\textup{(A1)}}\quad a_i<b_i=0;\qquad \mbox{\textup{(A2)}}\quad b_i>0;
\]
\item[(B)] for each $i \in\{1,2,\ldots,r\}$, $q_i+2a_i>q-t$;
\item[(C)] $N+2a_e > p-2 \sum_{i=1}^r a_i I_{(-\infty,0)}(a_i)$;
\item[(D)] $2b_e + \mathrm{SSE} > 0$.
\end{longlist}
\end{theorem}

If $m(y)<\infty$, then the posterior density is well defined
(i.e., proper) and is given by $\pi(\theta, \sigma^2|y) =
\pi^*(\theta,\sigma^2|y)  / m(y)$, but it is intractable in the
sense that posterior expectations cannot be computed in closed form,
nor even by classical Monte Carlo methods. However, there is a simple
two-step Gibbs sampler that can be used to approximate the intractable
posterior expectations. This Gibbs sampler simulates a Markov chain,
$\{(\theta_n,\sigma^2_n)\}_{n=0}^\infty$, that lives on $\X=
\mathbb{R}^{p+q} \times\mathbb{R}_+^{r+1}$, and has invariant density
$\pi(\theta, \sigma^2|y)$. If the current state of the chain is
$(\theta_n,\sigma^2_n)$, then the next state,
$(\theta_{n+1},\sigma^2_{n+1})$, is simulated using the usual two
steps. Indeed, we draw $\theta_{n+1}$ from
$\pi(\theta|\sigma^2_n,y)$, which is a $(p+q)$-dimensional
multivariate normal density, and then we draw $\sigma^2_{n+1}$ from
$\pi(\sigma^2|\theta_{n+1},y)$, which is a product of $r+1$ univariate
inverted gamma densities. The exact forms of these conditional
densities are given in Section~\ref{secsampler}.

Because the Gibbs--Markov chain is Harris ergodic (see
Section~\ref{secsampler}), we can use it to construct consistent
estimates of intractable posterior expectations. For $k>0$, let
$L_k(\pi)$ denote the set of functions $g\dvtx \mathbb{R}^{p+q} \times
\mathbb{R}_+^{r+1} \rightarrow\mathbb{R}$ such that
\[
E_\pi|g|^k:= \int_{\mathbb{R}_+^{r+1}} \int
_{\mathbb{R}^{p+q}} \bigl|g\bigl(\theta,\sigma^2\bigr)\bigr|^k
\pi\bigl(\theta, \sigma^2|y\bigr) \,d\theta \,d\sigma^2 <
\infty.
\]
If $g \in L_1(\pi)$, then the ergodic theorem implies that the average
\[
\overline{g}_m:= \frac{1}{m} \sum_{i=0}^{m-1}
g\bigl(\theta_i,\sigma_i^2\bigr)
\]
is a strongly consistent estimator of $E_\pi
g$, no matter how the chain is started. Of course, in practice, an
estimator is only useful if it is possible to compute an associated
(probabilistic) bound on the difference between the estimate and the
truth. Typically, this bound is based on a standard error. All
available methods of computing a valid asymptotic standard error for
$\overline{g}_m$ are based on the existence of a central limit theorem
(CLT) for $\overline{g}_m$ [see,
e.g., \citet
{joneharacaffneat2006,bednlatu2007,flegharajone2008,flegjone2010}].
Unfortunately, even if $g \in L_k(\pi)$ for all $k>0$, Harris
ergodicity is not enough to guarantee the existence of a CLT for
$\overline{g}_m$ [see, e.g., Roberts and Rosenthal
(\citeyear{roberose1998,roberose2004})].
The standard method of establishing the existence of CLTs is to prove
that the underlying Markov chain converges at a geometric rate.

Let ${\cal B}(\X)$ denote the Borel sets in $\X$, and let $P^n\dvtx  \X
\times{\cal B}(\X) \rightarrow[0,1]$ denote the $n$-step Markov
transition function of the Gibbs--Markov chain. That is,
$P^n ((\theta,\sigma^2),A  )$ is the probability that
$(\theta_n,\sigma^2_n) \in A$, given that the chain is started at
$(\theta_0,\sigma^2_0) = (\theta,\sigma^2)$. Also, let $\Pi(\cdot)$
denote the posterior distribution. The chain is called geometrically
ergodic if there exist a function $M\dvtx  \X\rightarrow[0,\infty)$ and a
constant $\varrho\in[0,1)$ such that, for all $(\theta,\sigma^2)
\in
\X$ and all $n=0,1,\ldots\,$, we have
\[
\bigl\| P^n \bigl(\bigl(\theta,\sigma^2\bigr),\cdot\bigr) -
\Pi(\cdot) \bigr\|_{\mathrm{TV}} \le M\bigl(\theta,\sigma^2\bigr)
\varrho^n,
\]
where \mbox{$\| \cdot\|_{\mathrm{TV}}$} denotes the total variation
norm. The relationship between geometric convergence and CLTs is
simple: if the chain is geometrically ergodic and $E_\pi
|g|^{2+\delta} < \infty$ for some $\delta>0$, then there is a CLT for
$\overline{g}_m$. Our main result (Theorem~\ref{thmmain} in
Section~\ref{secmain}) provides conditions under which the Gibbs--Markov chain is geometrically ergodic. The conditions of
Theorem~\ref{thmmain} are not easy to interpret, and checking them
may require some nontrivial numerical work. On the other hand, the
following corollary to Theorem~\ref{thmmain} is a slightly weaker
result whose conditions are very easy to check and understand.
%
\begin{corollary}
\label{corgeoerg}
Assume that $\operatorname{rank}(X)=p$. If the following four conditions
hold, then the Gibbs--Markov chain is geometrically ergodic.
\begin{longlist}[(C$'$)]
\item[(A)] For each $i \in\{1,2,\ldots,r\}$, one of the following
holds:
\[
\mbox{\textup{(A1)}}\quad a_i<b_i=0;\qquad \mbox{\textup{(A2)}}\quad b_i>0;
\]
\item[(B$'$)] for each $i \in\{1,2,\ldots,r\}$, $q_i+2a_i>q-t+2$;
\item[(C$'$)] $N+2a_e>p+t+2$;
\item[(D)] $2b_e + \mathrm{SSE} > 0$.
\end{longlist}
\end{corollary}

As we explain in Section~\ref{secsampler}, the best result
we could possibly hope to obtain is that the Gibbs--Markov chain is
geometrically ergodic whenever the posterior is proper. With this in
mind, note that the conditions of Corollary~\ref{corgeoerg} are very
close to the conditions for propriety given in
Theorem~\ref{thmpropriety}. In fact, the former imply the latter.
To see this, assume that (A), (B$'$), (C$'$) and (D) all hold.
Then, obviously, (B) holds, and all that remains is to show that
(C) holds. This would follow immediately if we could establish that
%
\begin{equation}
\label{eqmm} t \ge- 2 \sum_{i=1}^r
a_i I_{(-\infty,0)}(a_i).
\end{equation}
We consider two cases. First, if $\sum_{i=1}^r I_{(-\infty,0)}(a_i) =
0$, then it follows that $-2 \sum_{i=1}^r a_i I_{(-\infty,0)}(a_i)=0$,
and (\ref{eqmm}) holds (since $t$ is nonnegative). On the other
hand, if $\sum_{i=1}^r I_{(-\infty,0)}(a_i) > 0$, then there is at
least one negative $a_i$, and (B$'$) implies that
\[
\sum_{i=1}^r (q_i+2a_i)
I_{(-\infty,0)}(a_i) > q-t.
\]
This inequality combined with the fact that $q = q_1 + \cdots+ q_r$
yields
\[
t > q - \sum_{i=1}^r (q_i+2a_i)
I_{(-\infty,0)}(a_i) \ge- 2 \sum_{i=1}^r
a_i I_{(-\infty,0)}(a_i),
\]
so (\ref{eqmm}) holds, and this completes the argument.

The strong similarity between the conditions of
Corollary~\ref{corgeoerg} and those of Theorem~\ref{thmpropriety}
might lead the reader to believe that the proofs of our results rely
somehow on Theorem~\ref{thmpropriety}. This is not the case,
however. In fact, we do not even assume posterior propriety before
embarking on our convergence rate analysis; see
Section~\ref{secmain}.\vadjust{\goodbreak}

The only other existing result on geometric convergence of Gibbs
samplers for linear mixed models with \textit{improper} priors is that
of \citet{tanhobe2009}, who considered (a slightly reparameterized
version of) the one-way random effects model (\ref{eqowm}) and priors
with $b = (b_e, b_1) = (0, 0)$. We show in Section~\ref{secowm} that
our Theorem~\ref{thmmain} (specialized to the one-way model) improves
upon the result of \citet{tanhobe2009} in the sense that our
sufficient conditions for geometric convergence are weaker. Moreover,
it is known in this case exactly which priors lead to proper
posteriors (when $\mathrm{SSE}>0$), and we use this fact to show that
our results can be very close to the best possible. For example, if
the standard diffuse prior is used, then the posterior is proper if
and only if $c \ge3$. On the other hand, our results imply that the
Gibbs--Markov chain is geometrically ergodic as long as $c \ge3$, and
the total sample size, $N=n_1+n_2+\cdots+n_c$, is at least $c+2$. The
extra condition that $N \ge c+2$ is extremely weak. Indeed,
$\mathrm{SSE}>0$ implies that $N \ge c+1$, so, for fixed $c \ge3$, our
condition for geometric ergodicity fails only in the single case where
$N=c+1$.

An analogue of Corollary~\ref{corgeoerg} for the GLMM with
\textit{proper} priors can be found in \citet{johnjone2010}. In
contrast with our results, one of their sufficient conditions for
geometric convergence is that $X^TZ=0$, which rarely holds in
practice. Overall, the proper and improper cases are similar, in the
sense that geometric ergodicity is established via geometric drift
conditions in both cases. However, the drift conditions are quite
disparate, and the analysis required in the improper case is
substantially more demanding. Finally, we note that the linear models
considered by \citet{paparobe2008} are substantively different from
ours because these authors assume that the variance components are
known.

The remainder of this paper is organized as follows.
Section~\ref{secsampler} contains a formal definition of the
Gibbs--Markov chain. The main convergence result is stated and proven in
Section~\ref{secmain}, and an application involving the two-way
random effects model is given in Section~\ref{sectwm}. In
Section~\ref{secowm}, we consider the one-way random effects model
and compare our conditions for geometric convergence with those of
\citet{tanhobe2009}. Finally, Section~\ref{secdiscuss} concerns an
interesting technical issue related to the use of improper
priors.\vspace*{-2pt}

\section{The Gibbs sampler}
\label{secsampler}

In this section, we formally define the Markov chain\vspace*{1pt}
underlying the Gibbs sampler, and state some of its properties. Recall
that $\theta = (\beta^T\enskip u^T)^T$, $\sigma^2 =  (
\sigma^2_e\enskip \sigma^2_{u_1}\enskip \cdots\enskip \sigma^2_{u_r}
)^T$ and $\pi^*(\theta,\sigma^2|y)$ is the potentially
improper, unnormalized posterior density defined at (\ref{eqpistar}).
Suppose that
%
\begin{equation}
\label{eqtint} \int_{\mathbb{R}^{p+q}} \pi^*\bigl(\theta,
\sigma^2|y\bigr) \,d\theta< \infty
\end{equation}
for all $\sigma^2$ outside a set of measure zero in
$\mathbb{R}_+^{r+1}$, and that
%
\begin{equation}
\label{eqsint} \int_{\mathbb{R}_+^{r+1}} \pi^*\bigl(\theta,
\sigma^2|y\bigr) \,d\sigma^2 < \infty\vadjust{\goodbreak}
\end{equation}
for all $\theta$ outside a set of measure zero in $\mathbb{R}^{p+q}$.
These two integrability conditions are necessary, but not sufficient,
for posterior propriety. (Keep in mind that it is not known exactly
which priors yield proper posteriors.) When (\ref{eqtint}) and
(\ref{eqsint}) hold, we can define conditional densities as follows:
\[
\pi\bigl(\theta|\sigma^2,y\bigr) = \frac{\pi^*(\theta,\sigma^2|y)}{\int_{\mathbb{R}^{p+q}}
\pi^*(\theta,\sigma^2|y) \,d\theta}
\quad\mbox{and}\quad \pi\bigl(\sigma^2|\theta,y\bigr) =
\frac{\pi^*(\theta,\sigma^2|y)}{\int_{\mathbb{R}_+^{r+1}}
\pi^*(\theta,\sigma^2|y) \,d\sigma^2}.
\]
Clearly, when the posterior is proper, these conditionals are the
usual ones based on $\pi(\theta,\sigma^2|y)$. When the posterior is
improper, they are incompatible conditional densities; that is, there is
no (proper) joint density that generates them. In either case, we can
run the Gibbs sampler as usual by drawing alternately from the two
conditionals. However, as we explain below, if the posterior is
improper, then the resulting Markov chain cannot be geometrically
ergodic. Despite this fact, we do not restrict attention to the cases
where the sufficient conditions for propriety in
Theorem~\ref{thmpropriety} are satisfied. Indeed, we hope to close
the gap that currently exists between the necessary and sufficient
conditions for propriety by finding weaker conditions than those in
Theorem~\ref{thmpropriety} that imply geometric ergodicity (and hence
posterior propriety).

We now provide a set of conditions that guarantee that the
integrability conditions are satisfied. Define
\[
\tilde{s} = \min \{ q_1+2a_1, q_2+2a_2,
\ldots, q_r+2a_r, N+2a_e \}.
\]
The proof of the following result is straightforward and is left to
the reader.

\begin{proposition}
\label{propminimal}
The following four conditions are sufficient for (\ref{eqtint}) and
(\ref{eqsint}) to hold:
\begin{longlist}[(S4)]
\item[(S1)] $\operatorname{rank}(X) = p$;
\item[(S2)] $\min\{b_1,b_2,\ldots,b_r\} \ge0$;
\item[(S3)] $2b_e + \mathrm{SSE} > 0$;
\item[(S4)] $\tilde{s}>0$.
\end{longlist}
\end{proposition}

Note that $\mathrm{SSE} = \mathrm{SSE}(X,Z,y) = \llVert {y - W \hat
{\theta}}\rrVert^2$, where $W = (X\enskip Z)$ and $\hat{\theta} = (W^T
W)^{-}W^Ty$. Therefore, if condition (S3) holds, then for all
$\theta\in\mathbb{R}^{p+q}$,
\[
2b_e + \llVert {y-W\theta}\rrVert^2 = 2b_e
+ \llVert {y-W\hat {\theta}}\rrVert^2 + \llVert {W\theta-W\hat{
\theta}}\rrVert^2 \ge2b_e + \mathrm{SSE} > 0.
\]
Note also that if $N > p+q$, then $\mathrm{SSE}$ is strictly positive
with probability one under the data generating model.

Assume now that (S1)--(S4) hold so that the conditional densities
are well defined. Routine manipulation of $\pi^*(\theta,\sigma^2|y)$
shows that $\pi(\theta|\sigma^2,y)$ is a multivariate normal density
with mean vector
\[
m = \lleft[ %
\matrix{\bigl(X^TX\bigr)^{-1}X^T
\bigl( I-\bigl(\sigma^2_e\bigr)^{-1} Z
Q^{-1} Z^T (I-P_X) \bigr)y
\cr
\bigl(
\sigma^2_e\bigr)^{-1} Q^{-1}
Z^T (I-P_X) y} %
\rright]
\]
and covariance matrix
\[
V = \lleft[ %
\matrix{ \sigma^2_e
\bigl(X^TX\bigr)^{-1} + R Q^{-1}R^T &
-R Q^{-1}
\cr
-Q^{-1} R^T & Q^{-1} }
\rright],
\]
where $Q = (\sigma^2_e)^{-1} Z^T (I-P_X) Z + D^{-1}$ and
$R=(X^TX)^{-1} X^TZ$.

Things are a bit more complicated for $\pi(\sigma^2|\theta,y)$ due to
the possible existence of a bothersome set of measure zero. Define $A
=  \{ i \in\{1,2,\ldots,r \}\dvtx  b_i = 0  \}$. If $A$ is empty,
then $\pi(\sigma^2|\theta,y)$ is well defined for every $\theta\in
\mathbb{R}^{p+q}$, and it is the following product of $r+1$ inverted
gamma densities:
\begin{eqnarray*}
\pi\bigl( \sigma^2 |\theta,y\bigr)
&=& f_{\mathrm{IG}} \biggl(\sigma^2_e;
\frac{N}{2}+a_e, b_e+\frac{\|{y-W\theta}\|^2}{2} \biggr)
\\
&&{}\times\prod_{i=1}^r f_{\mathrm{IG}} \biggl(
\sigma^2_{u_i};\frac{q_i}{2}+a_i,
b_i+ \frac{\|{u_i}\|^2}{2} \biggr),
\end{eqnarray*}
where
\[
f_{\mathrm{IG}}(v;c,d) = \cases{\displaystyle \frac{d^c}{\Gamma(c) v^{c+1}} e^{-d/v}, &\quad $v > 0$,
\vspace*{2pt}\cr
0, &\quad $v \le0$, }
\]
for $c,d>0$. On the other hand, if $A$ is nonempty, then
\[
\int_{\mathbb{R}_+^{r+1}} \pi^*\bigl(\theta,\sigma^2|y\bigr) \,d
\sigma^2 = \infty,
\]
whenever $\theta\in{\cal N}:=  \{ \theta\in\mathbb{R}^{p+q}\dvtx
\prod_{i \in A} \llVert {u_i}\rrVert  = 0  \}$. The fact that
$\pi(\sigma^2|\theta,y)$ is not defined when $\theta\in{\cal N}$ is
irrelevant from a simulation standpoint because the probability of
observing a $\theta$ in ${\cal N}$ is zero. However, in order to
perform a theoretical analysis, the Markov transition density (Mtd) of
the Gibbs Markov chain must be defined for \textit{every} $\theta\in
\mathbb{R}^{p+q}$. Obviously, the Mtd can be defined arbitrarily on a
set of measure zero. Thus, for $\theta\notin{\cal N}$, we define
$\pi(\sigma^2|\theta,y)$ as in the case where $A$ is empty, while if
$\theta\in{\cal N}$, we define it to be $f_{\mathrm{IG}}
(\sigma^2_e;1,1) \prod_{i=1}^r
f_{\mathrm{IG}}(\sigma^2_{u_i};1,1)$. Note that this definition
can also be used when $A$ is empty if we simply define ${\cal N}$ to
be $\varnothing$ in that case.

The Mtd of the Gibbs--Markov chain,
$\{(\theta_n,\sigma_n^2)\}_{n=0}^\infty$, is defined as
\[
k\bigl(\theta,\sigma^2|\tilde{\theta},\tilde{\sigma}^2
\bigr) = \pi\bigl(\sigma^2|\theta,y\bigr) \pi\bigl(\theta|\tilde{
\sigma}^2,y\bigr).
\]
It is easy to see that the chain is $\psi$-irreducible, and that
$\pi^*(\theta,\sigma^2|y)$ is an invariant density. It follows that
the chain is positive recurrent if and only if the posterior is proper
[\citet{meyntwee1993}, Chapter 10]. Since a geometrically ergodic
chain is necessarily positive recurrent, the Gibbs--Markov chain cannot
be geometrically ergodic when the posterior is improper. The point
here is that conditions implying geometric ergodicity also imply
posterior propriety.

The marginal sequences, $\{\theta_n\}_{n=0}^\infty$ and
$\{\sigma_n^2\}_{n=0}^\infty$, are themselves Markov chains;
see, for example, \citet{liuwongkong1994}. The $\sigma^2$-chain lives
on $\mathbb{R}_+^{r+1}$ and has Mtd given by
\[
k_1\bigl(\sigma^2|\tilde{\sigma}^2\bigr) =
\int_{\mathbb{R}^{p+q}} \pi\bigl(\sigma^2|\theta,y\bigr) \pi
\bigl(\theta|\tilde{\sigma}^2,y\bigr) \,d\theta
\]
and invariant density $\int_{\mathbb{R}^{p+q}}
\pi^*(\theta,\sigma^2|y) \,d\theta$. Similarly, the $\theta$-chain
lives on $\mathbb{R}^{p+q}$ and has Mtd
\[
k_2(\theta|\tilde{\theta}) = \int_{\mathbb{R}_+^{r+1}} \pi\bigl(
\theta|\sigma^2,y\bigr) \pi\bigl(\sigma^2|\tilde{\theta},y
\bigr) \,d\sigma^2
\]
and invariant density $\int_{\mathbb{R}_+^{r+1}}
\pi^*(\theta,\sigma^2|y) \,d\sigma^2$. Since the two marginal chains
are also $\psi$-irreducible, they are positive recurrent if and only
if the posterior is proper. Moreover, when the posterior is proper,
routine calculations show that all three chains are Harris ergodic;
that is, $\psi$-irreducible, aperiodic and positive Harris recurrent;
see \citet{roma2012} for details. An important fact that we will
exploit is that geometric ergodicity is a solidarity property for the
three chains $\{(\theta_n,\sigma_n^2)\}_{n=0}^\infty$,
$\{\theta_n\}_{n=0}^\infty$ and $\{\sigma_n^2\}_{n=0}^\infty$; that
is, either all three are geometric or none of them is
[\citet{liuwongkong1994,roberose2001,diackharsalo2008}]. In the
next section, we prove that the Gibbs--Markov chain converges at a
geometric rate by proving that one of the marginal chains does.

\section{The main result}
\label{secmain}

In order to state the main result, we need a bit more notation.
For $i \in\{1,\ldots,r\}$, define $R_i$ to be the $q_i
\times q$ matrix of 0's and 1's such that $R_i u = u_i$. In other
words, $R_i$ is the matrix that \textit{extracts} $u_i$ from $u$.
Here is our main result.
%
\begin{theorem}
\label{thmmain}
Assume that \textup{(S1)--(S4)} hold so that the Gibbs sampler is well
defined. If the following two conditions hold, then the Gibbs--Markov
chain is geometrically ergodic:
\begin{longlist}[(2)]
\item[(1)] For each $i \in\{1,2,\ldots,r\}$, one of the following holds:
\[
\mbox{\textup{(i)}}\quad a_i<b_i=0;\qquad \mbox{\textup{(ii)}}\quad b_i>0.
\]
\item[(2)] There exists an $s \in(0,1] \cap(0,\tilde{s}/2)$ such that
%
\begin{equation}
\label{eqcondition2a} 2^{-s} (p+t)^s \frac{\Gamma ( {N}/{2}+a_e-s  )}{\Gamma
( {N}/{2}+a_e  )} < 1
\end{equation}
and
%
\begin{equation}
\label{eqcondition2b} 2^{-s} \sum_{i=1}^r
\biggl\{\frac{\Gamma ( {q_i}/{2}+a_i-s
)}{\Gamma ( {q_i}/{2}+a_i  )} \biggr\} \bigl( \tr \bigl( R_i
(I-P_ {Z^T (I-P_X) Z}) R_i^T \bigr) \bigr)^s <
1,
\end{equation}
where $t = \operatorname{rank}  ( Z^T (I-P_X) Z  )$ and $P_{Z^T (I-P_X)
Z}$ is the projection onto the column space of $Z^T (I-P_X) Z$.\vadjust{\goodbreak}
\end{longlist}
\end{theorem}
%
\begin{remark}
It is important to reiterate that, by themselves, (S1)--(S4)
\textit{do not} imply that the posterior density is proper. Of
course, if conditions (1) and (2) in Theorem~\ref{thmmain} hold as
well, then the chain is geometric, so the posterior is necessarily
proper.
\end{remark}
%
\begin{remark}
A numerical search could be employed to check the second condition
of Theorem~\ref{thmmain}. Indeed, one could evaluate the left-hand
sides of (\ref{eqcondition2a}) and (\ref{eqcondition2b}) at all
values of $s$ on a fine grid in the interval $(0,1] \cap
(0,\tilde{s}/2)$. The goal, of course, would be to find a single
value of $s$ at which both (\ref{eqcondition2a}) and
(\ref{eqcondition2b}) are satisfied. It can be shown that, if
there does exist an $s \in(0,1] \cap(0,\tilde{s}/2)$ such that
(\ref{eqcondition2a}) and (\ref{eqcondition2b}) hold, then
$N+2a_e>p+t$ and, for each $i=1,2,\ldots,r$, $q_i+2a_i > \tr(R_i
(I-P_ {Z^T (I-P_X) Z}) R_i^T)$. Thus, it would behoove the user to
verify these simple conditions before engaging in any numerical
work.
\end{remark}
%
\begin{remark}
When evaluating (\ref{eqcondition2b}), it may be helpful to write
$P_ {Z^T (I-P_X) Z}$ as $U^T P_\Lambda U$, where $U$ and $\Lambda$
are the orthogonal and diagonal matrices, respectively, in the
spectral decomposition of $Z^T(I-P_X)Z$. That is, $U$ is a
$q$-dimensional orthogonal matrix and $\Lambda$ is a diagonal matrix
containing the eigenvalues of $Z^T (I-P_X) Z$. Of course, the
projection $P_\Lambda$ is a $q \times q$ binary diagonal matrix
whose $i$th diagonal element is 1 if and only if the $i$th diagonal
element of $\Lambda$ is positive.
\end{remark}
%
\begin{remark}
\label{remre1}
Note that
\begin{eqnarray*}
\sum_{i=1}^r \tr \bigl( R_i
(I-P_ {Z^T (I-P_X) Z}) R_i^T \bigr) & = & \tr \Biggl[
(I-P_ {Z^T (I-P_X) Z}) \Biggl( \sum_{i=1}^r
R_i^T R_i \Biggr) \Biggr]
\\
& = & \tr(I-P_ {Z^T (I-P_X) Z})
\\
& = & \rank(I-P_
{Z^T (I-P_X) Z})
\\
& = & q-t.
\end{eqnarray*}
Moreover, when $r>1$, the matrix $I-P_{Z^T (I-P_X) Z}$ has
$q=q_1+q_2+\cdots+q_r$ diagonal elements, and the (nonnegative)
term $\tr(R_i (I-P_ {Z^T (I-P_X) Z}) R_i^T)$ is simply the sum of
the $q_i$ diagonal elements that correspond to the $i$th random
factor.
\end{remark}
%
\begin{remark}
Recall from the \hyperref[secintro]{Introduction} that Corollary~\ref{corgeoerg}
provides an alternative set of sufficient conditions for geometric
ergodicity that are harder to satisfy, but easier to check. A proof
of Corollary~\ref{corgeoerg} is given at the end of this section.
\end{remark}

We will prove Theorem~\ref{thmmain} indirectly by proving that the
$\sigma^2$-chain is geometrically ergodic (when the conditions of
Theorem~\ref{thmmain} hold). This is accomplished by establishing a
\textit{geometric drift condition} for the $\sigma^2$-chain.
%
\begin{proposition}
\label{propgd}
Assume that \textup{(S1)--(S4)} hold so that the Gibbs sampler is well
defined. Under the two conditions of Theorem~\ref{thmmain}, there
exist a $\rho\in[0,1)$ and a finite constant $L$ such that, for
every $\tilde{\sigma}^2 \in\mathbb{R}_+^{r+1}$,
%
\begin{equation}
\label{eqdrift} E \bigl( v \bigl(\sigma^2\bigr) | \tilde{
\sigma}^2 \bigr) \le\rho v \bigl(\tilde{\sigma}^2\bigr) +
L,
\end{equation}
where the drift function is defined as
\[
v \bigl(\sigma^2\bigr) = \alpha\bigl(\sigma_e^2
\bigr)^s + \sum_{i=1}^r \bigl(
\sigma_{u_i}^2\bigr)^s + \alpha\bigl(
\sigma_e^2\bigr)^{-c} + \sum
_{i=1}^r \bigl( \sigma_{u_i}^2
\bigr)^{-c},
\]
and $\alpha$ and $c$ are positive constants. Hence, under the two
conditions of Theorem~\ref{thmmain}, the $\sigma^2$-chain is
geometrically ergodic.
\end{proposition}
%
\begin{remark}
The formulas for $\rho=\rho(\alpha,s,c)$ and $L=L(\alpha,s,c)$ are
provided in the proof, as is a set of acceptable values for the pair
$(\alpha,c)$. Recall that the value of $s$ is given to us in the
hypothesis of Theorem~\ref{thmmain}.
\end{remark}
\begin{pf*}{Proof of Proposition~\ref{propgd}}
The proof has two parts. In part I, we establish the validity of
the geometric drift condition, (\ref{eqdrift}). In part II, we use
results from \citet{meyntwee1993} to show that (\ref{eqdrift})
implies geometric ergodicity of the $\sigma^2$-chain.

Part I. By conditioning on $\theta$ and iterating, we can express
$E(v(\sigma^2)|\tilde{\sigma}^2)$ as
\[
E \Biggl[ \alpha E \bigl( \bigl(\sigma^2_e
\bigr)^s | \theta \bigr) + E \Biggl( \sum
_{i=1}^r \bigl(\sigma^2_{u_i}
\bigr)^s | \theta \Biggr) + \alpha E \bigl( \bigl(
\sigma^2_e\bigr)^{-c} | \theta \bigr) + E
\Biggl( \sum_{i=1}^r \bigl(
\sigma^2_{u_i}\bigr)^{-c} | \theta \Biggr) \bigg|
\tilde{\sigma}^2 \Biggr].
\]
We now develop upper bounds for each of the four terms inside the
square brackets. Fix $s \in S:= (0,1] \cap(0,\tilde{s}/2)$, and
define
\[
G_0(s) = 2^{-s} \frac{\Gamma ( {N}/{2}+a_e-s  )}{\Gamma
( {N}/{2}+a_e  )}
\]
and, for each $i \in\{1,2,\ldots,r\}$, define
\[
G_i(s) = 2^{-s} \frac{\Gamma ( {q_i}/{2}+a_i-s
)}{\Gamma
( {q_i}/{2}+a_i  )}.
\]
Note that, since $s \in(0,1]$, $(x_1 + x_2)^s \le x_1^s + x_2^s$
whenever $x_1,x_2 \ge0$. Thus,
\begin{eqnarray*}
E \bigl( \bigl(\sigma^2_e\bigr)^s | \theta
\bigr) & = & 2^s G_0(s) \biggl( b_e +
\frac{\llVert {y-W\theta}\rrVert^2}{2} \biggr)^s
\\
& \le & 2^s G_0(s) \biggl[ b_e^s
+ \biggl( \frac{\llVert {y-W\theta}\rrVert^2}{2} \biggr)^s \biggr]
\\
& = & G_0(s) \bigl( \llVert {y-W\theta}\rrVert^2
\bigr)^s + 2^s G_0(s) b_e^s.
\end{eqnarray*}
Similarly,
\[
E \bigl( \bigl(\sigma^2_{u_i}\bigr)^s | \theta
\bigr) = 2^s G_i(s) \biggl( b_i +
\frac{\llVert {u_i}\rrVert^2}{2} \biggr)^s \le G_i(s) \bigl( \llVert
{u_i}\rrVert^2 \bigr)^s + 2^s
G_i(s) b_i^s.\vadjust{\goodbreak}
\]
Now, for any $c>0$, we have
\begin{eqnarray*}
E \bigl( \bigl(\sigma^2_e\bigr)^{-c} | \theta
\bigr) & = & 2^{-c} G_0(-c) \biggl( b_e +
\frac{\llVert {y-W\theta}\rrVert^2}{2} \biggr)^{-c}
\\
& \le & 2^{-c} G_0(-c) \biggl( b_e +
\frac{\mathrm{SSE}}{2} \biggr)^{-c}
\end{eqnarray*}
and, for each $i \in\{1,2,\ldots,r\}$,
\begin{eqnarray*}
E \bigl( \bigl(\sigma^2_{u_i}\bigr)^{-c} | \theta
\bigr) & = & 2^{-c} G_i(-c) \biggl( b_i +
\frac{\llVert {u_i}\rrVert^2}{2} \biggr)^{-c}
\\
& \le & G_i(-c) \bigl[ \bigl( \llVert {u_i}
\rrVert^2 \bigr)^{-c} I_{\{0\}
}(b_i) +
(2b_i)^{-c} I_{(0,\infty)}(b_i) \bigr].
\end{eqnarray*}
Recall that $A = \{i\dvtx  b_i=0 \}$, and note that $E  ( \sum_{i=1}^r
(\sigma^2_{u_i})^{-c}  | \theta )$ can be bounded above by a
constant if $A$ is empty. Thus, we consider the case in which $A$ is
empty separately from the case where $A \ne\varnothing$. We begin with
the latter, which is the more difficult case.

\textit{Case} I: $A$ is nonempty. Combining the four bounds above (and
applying Jensen's inequality twice), we have
%
\begin{eqnarray}
\label{eqpartial}\quad
E\bigl( v\bigl(\sigma^2\bigr)| \tilde{
\sigma}^2\bigr) &\le&\alpha G_0(s) \bigl[ E \bigl(
\llVert {y-W\theta}\rrVert^2 | \tilde{\sigma}^2 \bigr)
\bigr]^s + \sum_{i=1}^r
G_i(s) \bigl[ E \bigl( \llVert {u_i}
\rrVert^2 | \tilde{\sigma}^2 \bigr) \bigr]^s
\nonumber\\[-8pt]\\[-8pt]
&&{} + \sum_{i \in A} G_i(-c) E \bigl[ \llVert
{u_i}\rrVert^{-2c} | \tilde{\sigma}^2 \bigr] +
\kappa(\alpha,s,c),
\nonumber
\end{eqnarray}
where
\begin{eqnarray*}
\kappa(\alpha,s,c) &=& \alpha2^s G_0(s)
b_e^s + 2^s \sum
_{i=1}^r G_i(s) b_i^s
+ \alpha2^{-c} G_0(-c) \biggl( b_e +
\frac{\mathrm{SSE}}{2} \biggr)^{-c}
\\
&&{} + \sum_{i\dvtx b_i>0} G_i(-c)
(2b_i)^{-c}.
\nonumber
\end{eqnarray*}
Appendix~\ref{appub1} contains a proof of the following inequality:
%
\begin{equation}
\label{eqfirstexpec}\quad E \bigl[ \llVert {y-W\theta}\rrVert^2 | \tilde{
\sigma}^2 \bigr] \le (p+t) \tilde{\sigma}_e^{2}
+ \bigl( \bigl\llVert {(I-P_X) y}\bigr\rrVert + \bigl\llVert
{(I-P_X) Z}\bigr\rrVert K \bigr)^2,
\end{equation}
where \mbox{$\llVert \cdot\rrVert $} with a matrix argument denotes
the Frobenius
norm, and the constant $K=K(X,Z,y)$ is defined and shown to be finite
in Appendix~\ref{appprelim}. It follows immediately that
\[
\bigl[ E \bigl( \llVert {y-W\theta}\rrVert^2 | \tilde {
\sigma}^2 \bigr) \bigr]^s \le(p+t)^s \bigl(
\tilde{\sigma}^2_e\bigr)^s + \bigl( \bigl
\llVert {(I-P_X) y}\bigr\rrVert + \bigl\llVert {(I-P_X)
Z}\bigr\rrVert K \bigr)^{2s}.
\]
In Appendix~\ref{appub2}, it is shown that, for each $i \in
\{1,2,\ldots,r\}$, we have
\[
E \bigl[ \llVert {u_i}\rrVert^2 | \tilde{
\sigma}^2 \bigr] \le\xi_i \tilde{\sigma}_e^2
+ \zeta_i \sum_{j=1}^r
\tilde{\sigma}^2_{u_j} + \bigl( \llVert {R_i}
\rrVert K \bigr)^2,
\]
where $\xi_i = \tr ( R_i (Z^T(I-P_X)Z)^+ R_i^T  )$, $\zeta_i =
\tr ( R_i (I-P_{Z^T(I-P_X)Z}) R_i^T  )$ and $A^+$ denotes the
Moore--Penrose inverse of the matrix $A$. It follows that
%
\begin{equation}
\bigl[ E \bigl( \llVert {u_i}\rrVert^2 | \tilde{
\sigma}^2 \bigr) \bigr]^s \le\xi^s_i
\bigl(\tilde{\sigma}^2_e\bigr)^s +
\zeta^s_i \sum_{j=1}^r
\bigl(\tilde{\sigma}^2_{u_j}\bigr)^s + \bigl(
\llVert {R_i}\rrVert K \bigr)^{2s}.
\end{equation}
In Appendix~\ref{appub3}, it is established that, for any $c \in
(0,1/2)$, and for each $i \in\{1,2,\ldots,r\}$, we have
%
\begin{equation}
\label{eqfirstexpinv} E \bigl[ \llVert {{u_i}}\rrVert^{-2c}
| \tilde{\sigma}^2 \bigr] \le2^{-c} \frac{\Gamma ( {q_i}/{2}-c  )}{\Gamma (
{q_i}/{2}  )}
\bigl[ \lambda^c_{\max} \bigl(\tilde{\sigma}^2_e
\bigr)^{-c} + \bigl(\tilde{\sigma}^2_{u_i}
\bigr)^{-c} \bigr],
\end{equation}
where $\lambda_{\max}$ denotes the largest eigenvalue of $Z^T(I-P_X)
Z$. Using (\ref{eqfirstexpec})--(\ref{eqfirstexpinv}) in
(\ref{eqpartial}), we have
%
\begin{eqnarray}
\label{eqpartial2} E\bigl(v\bigl(\sigma^2\bigr)| \tilde{
\sigma}^2\bigr) &\le&\alpha \biggl( \delta_1(s) +
\frac{\delta_2(s)}{\alpha} \biggr) \bigl(\tilde{\sigma}^2_e
\bigr)^s + \delta_3(s) \sum_{j=1}^r
\bigl(\tilde{\sigma}^2_{u_j}\bigr)^s
\nonumber\\[-8pt]\\[-8pt]
&&{} + \alpha \frac{\delta_4(c)}{\alpha} \bigl(\tilde{\sigma}^2_e
\bigr)^{-c} + \delta_5(c) \sum
_{j \in A} \bigl(\tilde{\sigma}^2_{u_j}
\bigr)^{-c} + L(\alpha,s,c),
\nonumber
\end{eqnarray}
where
\begin{eqnarray*}
\delta_1(s)&:=& G_0(s) (p+t)^s,\qquad
\delta_2(s):= \sum_{i=1}^r
\xi_i^s G_i(s),\qquad
\delta_3(s):= \sum_{i=1}^r
\zeta_i^s G_i(s),
\\
\delta_4(c)&:=& 2^{-c} \lambda^c_{\max}
\sum_{i \in A} G_i(-c) \frac{\Gamma({q_i}/{2}-c)}{\Gamma({q_i}/{2})},
\\
\delta_5(c)&:=& 2^{-c} \max_{i \in A} \biggl[
G_i(-c) \frac{\Gamma({q_i}/{2}-c)}{\Gamma({q_i}/{2})} \biggr]
\end{eqnarray*}
and
\begin{eqnarray*}
L(\alpha,s,c) &=& \kappa(\alpha,s,c)
+ \alpha G_0(s) \bigl( \bigl\llVert {(I-P_X) y}\bigr
\rrVert + \bigl\llVert {(I-P_X) Z}\bigr\rrVert K
\bigr)^{2s}\\
&&{}  + \sum_{i=1}^rG_i(s)
\bigl( \llVert {R_i}\rrVert K \bigr)^{2s}.
\end{eqnarray*}
Hence,
\[
E\bigl(v\bigl(\sigma^2\bigr)| \tilde{\sigma}^2\bigr) \le
\rho(\alpha,s,c) v\bigl(\tilde{\sigma}^2\bigr) + L(\alpha,s,c),
\]
where
\[
\rho(\alpha,s,c) = \max \biggl\{ \delta_1(s) + \frac{\delta_2(s)}{\alpha},
\delta_3(s), \frac{\delta_4(c)}{\alpha}, \delta_5(c) \biggr\}.
\]
We must now show that there exists a triple $(\alpha,s,c) \in
\mathbb{R}_+ \times S \times(0,1/2)$ such that $\rho(\alpha,s,c)<1$.
We begin by demonstrating that, if $c$ is small enough, then
$\delta_5(c) < 1$. Define $\tilde{a} = -\max_{i \in A} a_i$. Also,
set $C = (0,1/2) \cap(0, \tilde{a})$. Fix $c \in C$ and note that
\[
\delta_5(c) = \max_{i \in A} \biggl[ \frac{\Gamma({q_i}/{2}+a_i+c)}
{\Gamma({q_i}/{2}+a_i)}
\frac{\Gamma({q_i}/{2}-c)}{\Gamma({q_i}/{2})} \biggr].
\]
For any $i \in A$, $c + a_i < 0$, and since $\tilde{s}>0$, it follows
that
\[
0 < \frac{q_i}{2} + a_i < \frac{q_i}{2} + a_i
+ c < \frac{q_i}{2}.
\]
But, $\Gamma(x-z)/\Gamma(x)$ is decreasing in $x$ for $x > z > 0$, so
we have
\[
\frac{\Gamma({q_i}/{2}+a_i)}{\Gamma({q_i}/{2}+a_i+c)} =
\frac{\Gamma({q_i}/{2}+a_i+c-c)}{\Gamma({q_i}/{2}+a_i+c)} >
\frac{\Gamma({q_i}/{2}-c)}{\Gamma({q_i}/{2})},
\]
and it follows immediately that $\delta_5(c) < 1$ whenever $c \in C$.
The two conditions of Theorem~\ref{thmmain} imply that there exists
an $s^\star\in S$ such that $\delta_1(s^\star) < 1$ and
$\delta_3(s^\star) < 1$. Let $c^\star$ be any point in $C$, and
choose $\alpha^\star$ to be any number larger than
\[
\max \biggl\{ \frac{\delta_2(s^\star)}{1-\delta_1(s^\star)}, \delta_4\bigl(c^\star
\bigr) \biggr\}.
\]
A simple calculation shows that $\rho(\alpha^\star,s^\star,c^\star
) <
1$, and this completes the argument for case I.

\textit{Case} II: $A = \varnothing$. Since we no longer have to deal with
$E ( \sum_{i=1}^r(\sigma^2_{u_i})^{-c}  | \theta )$, bound
(\ref{eqpartial2}) becomes
\[
E\bigl(v\bigl(\sigma^2\bigr)| \tilde{\sigma}^2\bigr) \le
\alpha \biggl( \delta_1(s) + \frac{\delta_2(s)}{\alpha} \biggr) \bigl(\tilde{
\sigma}^2_e\bigr)^s + \delta_3(s)
\sum_{j=1}^r \bigl(\tilde{
\sigma}^2_{u_j}\bigr)^s + L(\alpha,s,c),
\]
and there is no restriction on $c$ other than $c>0$. [Note that the
constant term $L(\alpha,s,c)$ requires no alteration when we move from
case I to case II.] Hence,
\[
E\bigl(v\bigl(\sigma^2\bigr)| \tilde{\sigma}^2\bigr) \le
\rho(\alpha,s) v\bigl(\tilde{\sigma}^2\bigr) + L(\alpha,s,c),
\]
where
\[
\rho(\alpha,s) = \max \biggl\{ \delta_1(s) + \frac{\delta_2(s)}{\alpha},
\delta_3(s) \biggr\}.
\]
We must now show that there exists a $(\alpha,s) \in\mathbb{R}_+
\times S$ such that $\rho(\alpha,s)<1$. As in case I, the two
conditions of Theorem~\ref{thmmain} imply that there exists an
$s^\star\in S$ such that $\delta_1(s^\star) < 1$ and
$\delta_3(s^\star) < 1$. Let $\alpha^\star$ be any number larger than
\[
\frac{\delta_2(s^\star)}{1-\delta_1(s^\star)}.
\]
A simple calculation shows that $\rho(\alpha^\star,s^\star) < 1$, and
this completes the argument for case II. This completes part I of the
proof.\vadjust{\goodbreak}

Part II. We begin by establishing that the $\sigma^2$-chain satisfies
certain properties. Recall that its Mtd is given by
\[
k_1\bigl(\sigma^2|\tilde{\sigma}^2\bigr) =
\int_{\mathbb{R}^{p+q}} \pi\bigl(\sigma^2|\theta,y\bigr) \pi
\bigl(\theta|\tilde{\sigma}^2,y\bigr) \,d\theta.
\]
Note that $k_1$ is strictly positive on $\mathbb{R}^{r+1}_+ \times
\mathbb{R}^{r+1}_+$. It follows that the $\sigma^2$-chain is
$\psi$-irreducible and aperiodic, and that its maximal irreducibility
measure is equivalent to Lebesgue measure on $\mathbb{R}^{r+1}_+$; for
definitions, see \citet{meyntwee1993}, Chapters 4 and 5. Let
$P_1$ denote the Markov transition function of the $\sigma^2$-chain;
that is, for any $\tilde{\sigma}^2 \in\mathbb{R}^{r+1}_+$ and any
Borel set $A$,
\[
P_1\bigl(\tilde{\sigma}^2,A\bigr) = \int
_A k_1\bigl(\sigma^2|\tilde{
\sigma}^2\bigr) \,d\sigma^2.
\]
We now demonstrate that the $\sigma^2$-chain is a Feller chain; that
is, for each fixed open set $O$, $P_1(\cdot,O)$ is a lower
semi-continuous function on $\mathbb{R}^{r+1}_+$. Indeed, let
$\{\tilde{\sigma}_m^2\}_{m=1}^\infty$ be a sequence in
$\mathbb{R}^{r+1}_+$ that converges to $\tilde{\sigma}^2 \in
\mathbb{R}^{r+1}_+$. Then
\begin{eqnarray*}
\liminf_{m \rightarrow\infty} P_1\bigl(\tilde{\sigma}^2_m,O
\bigr) & = & \liminf_{m \rightarrow\infty} \int_O k_1
\bigl(\sigma^2|\tilde{\sigma}^2_m\bigr) \,d
\sigma^2
\\
& = &\liminf_{m \rightarrow\infty} \int_O \biggl[ \int
_{\mathbb{R}^{p+q}} \pi\bigl(\sigma^2|\theta,y\bigr) \pi\bigl(
\theta|\tilde{\sigma}^2_m,y\bigr) \,d\theta \biggr] \,d
\sigma^2
\\
& \ge &\int_O \int_{\mathbb{R}^{p+q}} \pi\bigl(
\sigma^2|\theta,y\bigr) \Bigl[ \liminf_{m \rightarrow\infty} \pi\bigl(\theta|
\tilde{\sigma}^2_m,y\bigr) \Bigr] \,d\theta \,d
\sigma^2
\\
& = &\int_O \biggl[ \int_{\mathbb{R}^{p+q}} \pi
\bigl(\sigma^2|\theta,y\bigr) \pi\bigl(\theta|\tilde{
\sigma}^2,y\bigr) \,d\theta \biggr] \,d\sigma^2
\\
& = & P_1\bigl(\tilde{\sigma}^2,O\bigr),
\end{eqnarray*}
where the inequality follows from Fatou's lemma, and the third
equality follows from the fact that $\pi(\theta|\sigma^2,y)$ is
continuous in $\sigma^2$; for a proof of continuity, see
\citet{roma2012}. We conclude that $P_1(\cdot,O)$ is lower
semi-continuous, so the $\sigma^2$-chain is Feller.

The last thing we must do before we can appeal to the results in
\citet{meyntwee1993} is to show that the drift function, $v(\cdot)$,
is unbounded off compact sets; that is, we must show that, for every
$d \in\mathbb{R}$, the set
\[
S_d = \bigl\{ \sigma^2 \in\mathbb{R}^{r+1}_+\dvtx
v\bigl(\sigma^2\bigr) \le d \bigr\}
\]
is compact. Let $d$ be such that $S_d$ is nonempty (otherwise $S_d$
is trivially compact), which means that $d$ and $d/\alpha$ must be
larger than 1. Since $v(\sigma^2)$ is a continuous function, $S_d$ is
closed in $\mathbb{R}^{r+1}_+$. Now consider the following set:
\[
T_d = \bigl[ (d/\alpha)^{-1/c}, (d/\alpha)^{1/s}
\bigr] \times \bigl[ d^{-1/c}, d^{1/s} \bigr] \times\cdots\times
\bigl[ d^{-1/c}, d^{1/s} \bigr].
\]
The set $T_d$ is compact in $\mathbb{R}_+^{r+1}$. Since $S_d \subset
T_d$, $S_d$ is a closed subset of a compact set in
$\mathbb{R}_+^{r+1}$, so it is compact in $\mathbb{R}_+^{r+1}$.
Hence, the drift function is unbounded off compact sets.

Since the $\sigma^2$-chain is Feller and its maximal irreducibility
measure is equivalent to Lebesgue measure on $\mathbb{R}^{r+1}_+$,
\citeauthor{meyntwee1993}'s (\citeyear{meyntwee1993}) Theorem
6.0.1 shows that every compact set in
$\mathbb{R}^{r+1}_+$ is petite. Hence, for each $d \in\mathbb{R}$,
the set $S_d$ is petite, so $v(\cdot)$ is unbounded off petite sets.
It now follows from the drift condition (\ref{eqdrift}) and an
application of \citeauthor{meyntwee1993}'s (\citeyear
{meyntwee1993}) Lemma 15.2.8 that condition
(iii) of \citeauthor{meyntwee1993}'s (\citeyear{meyntwee1993})
Theorem 15.0.1 is satisfied, so the
$\sigma^2$-chain is geometrically ergodic. This completes part II of
the proof.
\end{pf*}

We end this section with a proof of Corollary~\ref{corgeoerg}.
\begin{pf*}{Proof of Corollary~\ref{corgeoerg}}
It suffices to show that, together, conditions (B$'$) and (C$'$) of
Corollary~\ref{corgeoerg} imply the second condition of
Theorem~\ref{thmmain}. Clearly, (B$'$) and (C$'$) imply that
$\tilde{s}/2 > 1$, so $(0,1] \cap(0,\tilde{s}/2) = (0,1]$. Take
$s^\star= 1$. Condition (C$'$) implies
\[
2^{-s^\star} (p+t)^{s^\star} \frac{\Gamma (
{N}/{2}+a_e-s^\star )}{\Gamma ( {N}/{2}+a_e  )} = \frac{p+t}{N+2a_e-2} <
1.
\]
Now, we know from Remark~\ref{remre1} that $\sum_{i=1}^r \tr (
R_i(I-P_{Z^T(I-P_X)Z})R_i^T  ) = q-t$. Hence,
\begin{eqnarray*}
&&
2^{-s^\star}  \sum_{i=1}^r \biggl\{
\frac{\Gamma (
{q_i}/{2}+a_i-s^\star )}{\Gamma ( {q_i}/{2}+a_i
)} \biggr\} \bigl( \tr \bigl( R_i
(I-P_{Z^T(I-P_X)Z}) R_i^T \bigr) \bigr)^{s^\star}
\\
&&\qquad= \sum_{i=1}^r \frac{\tr ( R_i (I-P_{Z^T(I-P_X)Z}) R_i^T
)}{q_i + 2a_i -2}
\\
&&\qquad \le\frac{\sum_{i=1}^r \tr ( R_i
(I-P_{Z^T(I-P_X)Z})R_i^T  )}{\min_{j \in\{1,2,\ldots,r\}} \{q_j + 2a_j
-2\}}
\\
&&\qquad = \frac{q-t}{\min_{j \in\{1,2,\ldots,r\}} \{q_j + 2a_j -2\}} < 1,
\end{eqnarray*}
where the last inequality follows from condition (B$'$).
\end{pf*}

\section{An application of the main result}
\label{sectwm}

In this section, we illustrate the application of
Theorem~\ref{thmmain} using the two-way random effects model with one
observation per cell. The model equation is
\[
Y_{ij} = \beta+ \alpha_i + \gamma_j +
\varepsilon_{ij},
\]
where $i=1, 2, \ldots, m$, $j=1, 2, \ldots, n$, the $\alpha_i$'s are i.i.d.
$N(0,\sigma^2_{\alpha})$, the $\gamma_j$'s are i.i.d.
$N(0,\sigma^2_{\gamma})$ and the $\varepsilon_{ij}$'s are i.i.d. $N(0,
\sigma^2_e)$. The $\alpha_i$'s, $\gamma_j$'s and $\varepsilon_{ij}$'s are
all independent.

We begin by explaining how to put this model in GLMM (matrix) form.
There are a total of $N = m \times n$ observations and we arrange them
using the usual (lexicographical) ordering
\[
Y = \pmatrix{ Y_{11} \cdots Y_{1n} & Y_{21} \cdots
Y_{2n} &\cdots& Y_{m1} \cdots Y_{mn} }^T.
\]
Since $\beta$ is a univariate parameter common to all of the
observations, $p=1$ and $X$ is an $N \times1$ column vector of ones,
which we denote by $1_N$. There are $r=2$ random factors with $q_1=m$
and $q_2=n$, so $q=m+n$. Letting $\otimes$ denote the Kronecker product,
we can write the $Z$ matrix as $(Z_1\enskip Z_2)$, where $Z_1 = I_m
\otimes1_n$ and $Z_2 = 1_m \otimes I_n$. We assume throughout this
section that $\mathrm{SSE}>0$.

We now examine conditions (\ref{eqcondition2a}) and
(\ref{eqcondition2b}) of Theorem~\ref{thmmain} for this particular
model. A routine calculation shows that $Z^T(I-P_X)Z$ is a block
diagonal matrix given by
\[
n \biggl(I_m-\frac{1}{m}J_m \biggr) \oplus m
\biggl(I_n-\frac
{1}{n}J_n \biggr),
\]
where $J_{d}$ is a $d \times d$ matrix of ones (and $\oplus$ is the
direct sum operator). It follows immediately that
\[
t = \operatorname{rank} \bigl( Z^T (I-P_X) Z \bigr) =
\operatorname{rank} \biggl(I_m-\frac{1}{m}J_m
\biggr) + \operatorname{rank} \biggl(I_n-\frac{1}{n}J_n
\biggr) = m+n-2.
\]
Hence, (\ref{eqcondition2a}) becomes
\[
2^{-s} (m+n-1)^s \frac{\Gamma ( {mn}/{2}+a_e-s
)}{\Gamma
( {mn}/{2}+a_e  )} < 1.
\]
Now, it can be shown that
\[
I-P_{Z^T(I-P_X)Z} = \biggl(\frac{1}{m}J_m \biggr) \oplus
\biggl(\frac{1}{n}J_n \biggr).
\]
Hence, we have
\[
\tr\bigl(R_1 (I-P_{Z^T(I-P_X)Z})R_1^T
\bigr) = \tr \biggl(\frac{1}{m}J_m \biggr) = 1
\]
and
\[
\tr\bigl(R_2 (I-P_{Z^T(I-P_X)Z}) R_2^T
\bigr) = \tr \biggl(\frac
{1}{n}J_n \biggr) = 1.
\]
Therefore, (\ref{eqcondition2b}) reduces to
\[
2^{-s} \biggl\{\frac{\Gamma ( {m}/{2}+a_1-s
)}{\Gamma ( {m}/{2}+a_1  )} + \frac{\Gamma ( {n}/{2}+a_2-s
)}{\Gamma ( {n}/{2}+a_2  )} \biggr\} < 1.
\]

Now consider a concrete example in which $m=5$, $n=6$, and the prior
is
\[
\frac{I_{\mathbb{R}_+}(\sigma_e^2) I_{\mathbb{R}_+}(\sigma_{\alpha}^2)
I_{\mathbb{R}_+}(\sigma_{\gamma}^2)}{\sigma^2_e
\sqrt{\sigma^2_{\alpha} \sigma^2_{\gamma}}}.
\]
So, we are taking $b_e=b_1=b_2=a_e=0$ and $a_1=a_2=-1/2$.
Corollary~\ref{corgeoerg} implies that the Gibbs--Markov chain is
geometrically ergodic whenever $m,n \ge6$, but this result is not
applicable when $m=5$ and $n=6$. Hence, we turn to
Theorem~\ref{thmmain}. In this case, $\tilde{s} = 4$, so we need to
find an $s \in(0,1]$ such that
\begin{eqnarray*}
&&
2^{-s} \max \biggl\{ (10)^s \frac{\Gamma ( {30}/{2}+0-s
)}{\Gamma ( {30}/{2}+0  )},\\
&&\qquad\hspace*{21.1pt}
\frac{\Gamma (
{5}/{2}+(-{1}/{2})-s  )}{\Gamma (
{5}/{2}+(-{1}/{2})  )} + \frac{\Gamma (
{6}/{2}+(-{1}/{2})-s  )}{\Gamma (
{6}/{2}+(-{1}/{2})  )} \biggr\} < 1.
\end{eqnarray*}
The reader can check that, when $s= 0.9$, the left-hand side is
approximately $0.87$. Therefore, Theorem~\ref{thmmain} implies that
the Gibbs--Markov chain is geometrically ergodic in this case.

\section{Specializing to the one-way random effects model}
\label{secowm}

The only other existing results on geometric convergence of Gibbs
samplers for linear mixed models with \textit{improper} priors are
those of \citet{tanhobe2009} (hereafter, T\&H). These authors
considered the one-way random effects model, which is a simple, but
important special case of the GLMM given in (\ref{eqGLMM}). In this
section, we show that our results improve upon those of T\&H.

Recall that the one-way model is given by
%
\begin{equation}
\label{eqowm2} Y_{ij} = \beta+ \alpha_i +
e_{ij},
\end{equation}
where $i=1,\ldots,c$, $j=1,\ldots,n_i$, the $\alpha_i$'s are i.i.d.
$\mathrm{N}(0,\sigma^2_\alpha)$, and the $e_{ij}$'s, which are
independent of the $\alpha_i$'s, are i.i.d. $\mathrm{N}(0,\sigma_e^2)$.
It is easy to see that (\ref{eqowm2}) is a special case of the GLMM.
Obviously, there are a total of $N = n_1 + \cdots+ n_c$ observations,
and we arrange them in a column vector with the usual ordering as
follows:
\[
Y = \pmatrix{ Y_{11} \cdots Y_{1n_1} & Y_{21} \cdots
Y_{2n_2} &\cdots& Y_{c1} \cdots Y_{cn_c}
}^T.
\]
As in the two-way model of Section~\ref{sectwm}, $\beta$ is a
univariate parameter common to all of the observations, so $p=1$ and
$X = 1_N$. Here there is only one random factor (with $c$ levels), so
$r=1$, $q=q_1=c$ and $Z = \bigoplus_{i=1}^c 1_{n_i}$. Of course, in this
case, $\mathrm{SSE} = \sum_{i=1}^c \sum_{j=1}^{n_i} (y_{ij} -
\overline{y}_i)^2$, where $\overline{y}_i = n_i^{-1} \sum_{j=1}^{n_i}
y_{ij}$. We assume throughout this section that $\mathrm{SSE}>0$.

We note that T\&H actually considered a slightly different
parameterization of the one-way model. In their version, $\beta$ does
not appear in the model equation~(\ref{eqowm2}), but rather as the
mean of the $u_i$'s. In other words, T\&H used the ``centered''
parameterization, whereas here we are using the ``noncentered''
parameterization. \citet{roma2012} shows that, because $\beta$ and
$(\alpha_1 \cdots\alpha_c)^T$ are part of a single ``block'' in the
Gibbs sampler, the centered and noncentered versions of the Gibbs
sampler converge at exactly the same rate.

T\&H considered improper priors for
$(\beta,\sigma_e^2,\sigma_\alpha^2)$ that take the form
\[
\bigl(\sigma^2_e\bigr)^{-(a_e+1)} I_{\mathbb{R}_+}
\bigl(\sigma^2_e\bigr) \bigl(\sigma^2_\alpha
\bigr)^{-(a_1+1)} I_{\mathbb{R}_+}\bigl(\sigma^2_\alpha
\bigr),
\]
and they showed that the Gibbs sampler for the one-way model is
geometrically ergodic if $a_1<0$ and
%
\begin{eqnarray}
\label{eqTHc}
N+2a_e &\ge& c+3 \quad\mbox{and}\nonumber\\[-8pt]\\[-8pt]
c \min \Biggl\{ \Biggl( \sum_{i=1}^c
\frac{n_i}{n_i+1} \Biggr)^{-1}, \frac{n^*}{N} \Biggr\} &<& 2 \exp
\biggl\{ \Psi \biggl( \frac{c}{2}+a_1 \biggr) \biggr\},\nonumber
\end{eqnarray}
where $n^* = \max\{n_1,n_2,\ldots,n_c\}$ and $\Psi(x) = \frac{d}{dx}
\log ( \Gamma(x)  )$ is the digamma function.

We now consider the implications of Theorem~\ref{thmmain} in the case
of the one-way model. First, $t = \operatorname{rank}  ( Z^T
(I-P_X) Z
)=c-1$. Combining this fact with Remark~\ref{remre1}, it follows
that the two conditions of Theorem~\ref{thmmain} will hold if $a_1<0$,
and there exists an $s \in(0,1) \cap(0,a_1+\frac{c}{2}) \cap
(0,a_e+\frac{N}{2})$ such that
\[
2^{-s} \max \biggl\{ c^s \frac{\Gamma ( {N}/{2}+a_e-s
)}{\Gamma ( {N}/{2}+a_e  )},
\frac{\Gamma (
{c}/{2}+a_1-s  )}{\Gamma ( {c}/{2}+a_1  )} \biggr\} < 1.
\]
\citet{roma2012} shows that such an $s$ does indeed exist (so the
Gibbs chain is geometrically ergodic) when
%
\begin{equation}
\label{eqRHc} N+2a_e \ge c+2 \quad\mbox{and}\quad
1 < 2 \exp \biggl\{ \Psi \biggl( \frac{c}{2}+a_1 \biggr)
\biggr\}.
\end{equation}
Now, it's easy to show that
\[
1 \le c \min \Biggl\{ \Biggl( \sum_{i=1}^c
\frac{n_i}{n_i+1} \Biggr)^{-1}, \frac{n^*}{N} \Biggr\}.
\]
Consequently, if (\ref{eqTHc}) holds, then so does (\ref{eqRHc}).
In other words, our sufficient conditions are weaker than those of
T\&H, so our result improves upon theirs. Moreover, in contrast with
the conditions of T\&H, our conditions do not directly involve the
group sample sizes, $n_1, n_2, \ldots, n_c$.

Of course, the best result possible would be that the Gibbs--Markov
chain is geometrically ergodic whenever the posterior is proper. Our
result is very close to the best possible in the important case where
the standard diffuse prior is used; that is, when $a_1 = -1/2$ and
$a_e = 0$. The posterior is proper in this case if and only if $c \ge
3$ [\citet{SunTsutHe2001}]. It follows from (\ref{eqRHc}) that the
Gibbs--Markov chain is geometrically ergodic as long as $c \ge3$, and
the total sample size, $N=n_1+n_2+\cdots+n_c$, is at least $c+2$.
This additional sample size condition is \textit{extremely weak}.
Indeed, the positivity of $\mathrm{SSE}$ implies that $N \ge c+1$, so,
for fixed $c \ge3$, our condition for geometric ergodicity fails only
in the single case where $N=c+1$. Interestingly, in this case, the
conditions of Corollary~\ref{corgeoerg} reduce to $c \ge5$ and $N
\ge c+3$.

\section{Discussion}
\label{secdiscuss}

Our decision to work with the $\sigma^2$-chain rather than the
$\theta$-chain was based on an important technical difference between
the two chains that stems from the fact that $\pi(\sigma^2|\theta,y)$
is not continuous in $\theta$ for each fixed $\sigma^2$ (when the set
$A$ is nonempty). Indeed, recall that, for $\theta\notin{\cal N}$,
\begin{eqnarray*}
\pi\bigl( \sigma^2 |\theta,y\bigr)
&=& f_{\mathrm{IG}} \biggl(\sigma^2_e;
\frac{N}{2}+a_e, b_e+\frac{\llVert {y-W\theta
}\rrVert^2}{2} \biggr)\\
&&{}\times
\prod_{i=1}^r f_{\mathrm{IG}} \biggl(
\sigma^2_{u_i};\frac{q_i}{2}+a_i,
b_i+ \frac{\llVert {u_i}\rrVert^2}{2} \biggr),\nonumber
\end{eqnarray*}
but for $\theta\in{\cal N}$,
\[
\pi\bigl(\sigma^2|\theta,y\bigr) = f_{\mathrm{IG}} \bigl(
\sigma^2_e;1,1\bigr) \prod_{i=1}^r
f_{\mathrm{IG}}\bigl(\sigma^2_{u_i};1,1\bigr).
\]
Also, recall that the Mtd of the $\sigma^2$-chain is given by
\[
k_1\bigl(\sigma^2|\tilde{\sigma}^2\bigr) =
\int_{\mathbb{R}^{p+q}} \pi\bigl(\sigma^2|\theta,y\bigr) \pi
\bigl(\theta|\tilde{\sigma}^2,y\bigr) \,d\theta.
\]
Since the set ${\cal N}$ has measure zero, the ``arbitrary part'' of
$\pi(\sigma^2|\theta,y)$ washes out of~$k_1$. However, the same
cannot be said for the $\theta$-chain, whose Mtd is given by
\[
k_2(\theta|\tilde{\theta}) = \int_{\mathbb{R}_+^{r+1}} \pi\bigl(
\theta|\sigma^2,y\bigr) \pi\bigl(\sigma^2|\tilde{\theta},y
\bigr) \,d\sigma^2.
\]
This difference between $k_1$ and $k_2$ comes into play when we
attempt to apply certain ``topological'' results from Markov chain
theory, such as those in Chapter 6 of \citet{meyntwee1993}. In
particular, in our proof that the $\sigma^2$-chain is a Feller chain
(which was part of the proof of Proposition~\ref{propgd}), we used
the fact that $\pi(\theta|\sigma^2,y)$ is continuous in $\sigma^2$ for
each fixed $\theta$. Since $\pi(\sigma^2|\theta,y)$ is not
continuous, we cannot use the same argument to prove that the
$\theta$-chain is Feller. In fact, we suspect that the $\theta$-chain
is not Feller, and if this is true, it means that our method of proof
will not work for the $\theta$-chain.

It is possible to circumvent the problem described above by removing
the set ${\cal N}$ from the state space of the $\theta$-chain.\vadjust{\goodbreak} In
this case, we are no longer required to define
$\pi(\sigma^2|\theta,y)$ for $\theta\in{\cal N}$, and since
$\pi(\sigma^2|\theta,y)$ is continuous (for fixed $\sigma^2$) on
$\mathbb{R}^{p+q} \setminus{\cal N}$, the Feller argument for the
$\theta$-chain will go through. On the other hand, the new state
space has ``holes'' in it, and this could complicate the search for a
drift function that is unbounded off compact sets. For example,
consider a toy drift function given by $v(x) = x^2$. This function is
clearly unbounded off compact sets when the state space is
$\mathbb{R}$, but not when the state space is $\mathbb{R} \setminus
\{0\}$. The modified drift function $v^*(x) = x^2 + 1/x^2$ is
unbounded off compact sets for the ``holey'' state space.

T\&H overlooked a set of measure zero (similar to our ${\cal N}$), and
this oversight led to an error in the proof of their main result
(Proposition 3). However, \citet{roma2012} shows that T\&H's proof
can be repaired and that their result is correct as stated. The fix
involves deleting the offending null set from the state space, and
adding a term to the drift function.

\begin{appendix}\label{app}
\section*{Appendix: Upper bounds}

\subsection{Preliminary results}
\label{appprelim}

Here is our first result.
%
\begin{lemma}
\label{lemmatracebounds}
The following inequalities hold for all $\sigma^2 \in
{\mathbb{R}_+^{r+1}}$ and all $i \in\{1,2,\ldots,r\}$:
\begin{longlist}[(3)]
\item[(1)] $Q^{-1} \preceq(Z^T (I-P_X) Z)^+ \sigma_e^{2} +
(I-P_{Z^T(I-P_X)Z}) (
\sum_{j=1}^r \sigma^2_{u_j}  )$;
\item[(2)] $\tr ( (I-P_X) Z Q^{-1} Z^T (I-P_X)  ) \le
\operatorname{rank}(Z^T (I-P_X) Z) \sigma_e^{2}$;
\item[(3)] $(R_i Q^{-1} R^T_i)^{-1} \preceq ( (\sigma^2_e)^{-1}
\lambda_{\max} + (\sigma^2_{u_i})^{-1}  ) I_{q_i}$.
\end{longlist}
\end{lemma}
\begin{pf}
Recall from Section~\ref{secmain} that $U^T \Lambda U$ is the
spectral decomposition of $Z^T(I-P_X) Z$, and that $P_\Lambda$ is a binary
diagonal matrix whose $i$th diagonal element is 1 if and only if the
$i$th diagonal element of $\Lambda$ is positive. Let $\sigma^2_\bullet=
\sum_{j=1}^r \sigma^2_{u_j}$. Since $ ( \sigma^2_\bullet
)^{-1} I_q \preceq D^{-1}$, we have
\[
\bigl(\sigma_e^{2}\bigr)^{-1} Z^T
(I-P_X) Z + \bigl(\sigma^2_\bullet
\bigr)^{-1} I_q \preceq\bigl(\sigma^2_e
\bigr)^{-1} Z^T (I-P_X) Z + D^{-1},
\]
and this yields
%
\begin{eqnarray}
\label{eqPlamb} Q^{-1} & = & \bigl( \bigl(\sigma^2_e
\bigr)^{-1} Z^T (I-P_X) Z + D^{-1}
\bigr)^{-1}
\nonumber
\\
& \preceq &\bigl( \bigl(\sigma^2_e\bigr)^{-1}
Z^T (I-P_X) Z + \bigl(\sigma^2_\bullet
\bigr)^{-1} I_q \bigr)^{-1}
\\
& = &U^T \bigl( \Lambda\bigl(\sigma^2_e
\bigr)^{-1} + I_q \bigl(\sigma^2_\bullet
\bigr)^{-1} \bigr)^{-1} U.\nonumber
\end{eqnarray}
Now let $\Lambda^+$ be a diagonal matrix whose diagonal elements,
$\{\lambda_i^+\}_{i=1}^q$, are given by
\[
\lambda_i^+ = \cases{ \lambda_i^{-1}, &\quad $
\lambda_i \ne0$,
\cr
0, &\quad $\lambda_i = 0$. }
\]
Note that, for each $i \in\{1,2,\ldots,r\}$, we have
\[
\frac{1}{\lambda_i (\sigma^2_e)^{-1} + (\sigma^2_\bullet)^{-1}} \le \lambda_i^{+}
\sigma_e^{2} +I_{\{0\}}(\lambda_i)
\sigma^2_\bullet.\vadjust{\goodbreak}
\]
This shows that
\[
\bigl( \Lambda\bigl(\sigma^2_e\bigr)^{-1} +
I_q \bigl(\sigma^2_\bullet\bigr)^{-1}
\bigr)^{-1} \preceq\Lambda^+ \sigma^2_e +
(I-P_\Lambda) \sigma^2_\bullet.
\]
Together with (\ref{eqPlamb}), this leads to
\begin{eqnarray*}
Q^{-1} &\preceq& U^T \bigl( \Lambda\bigl(
\sigma^2_e\bigr)^{-1} + I_q \bigl(
\sigma^2_\bullet\bigr)^{-1} \bigr)^{-1}
U \preceq U^T \bigl( \Lambda^+ \sigma^2_e +
(I-P_\Lambda) \sigma^2_\bullet \bigr) U
\\
& = &\bigl(Z^T(I-P_X) Z\bigr)^+
\sigma^2_e + U^T (I-P_\Lambda) U
\sigma^2_\bullet.
\end{eqnarray*}
So to prove the first statement, it remains to show that $U^T
(I-P_\Lambda) U = I- P_{Z^T(I-P_X)Z}$. But\vspace*{1pt} notice that letting $A =
Z^T(I-P_X)Z$ and using its spectral decomposition, we have
\begin{eqnarray*}
A\bigl(A^TA\bigr)^+A^T &=& U^T\Lambda U
\bigl(U^T\Lambda^T\Lambda U\bigr)^+ U^T\Lambda
U \\
&=& U^T\Lambda\bigl(\Lambda^T\Lambda\bigr)^+ \Lambda U =
U^TP_\Lambda U,
\end{eqnarray*}
which implies that
\[
I- P_{Z^T(I-P_X)Z} =I - A\bigl(A^TA\bigr)^+A^T = I -
U^T P_\Lambda U = U^T (I-P_\Lambda) U.
\]
The proof of the first statement is now complete. Now let $\tilde{Z}
= (I-P_X) Z$. Multiplying the first statement on the left and the
right by $\tilde{Z}$ and $\tilde{Z}^T$, respectively, and then taking
traces yields
%
\begin{equation}
\label{eqn1} \tr \bigl( \tilde{Z} Q^{-1} \tilde{Z}^T
\bigr) \le\tr \bigl( \tilde{Z} \bigl(\tilde{Z}^T \tilde{Z}\bigr)^+
\tilde{Z}^T \bigr) \sigma_e^2 + \tr \bigl(
\tilde{Z} U^T (I-P_\Lambda) U \tilde{Z}^T \bigr)
\sigma^2_\bullet.
\end{equation}
Since $(\tilde{Z}^T \tilde{Z}) (\tilde{Z}^T \tilde{Z})^+$ is
idempotent, we have
\begin{eqnarray*}
\tr \bigl( \tilde{Z} \bigl(\tilde{Z}^T \tilde{Z}\bigr)^+
\tilde{Z}^T \bigr) &=& \tr \bigl( \tilde{Z}^T \tilde{Z}
\bigl(\tilde{Z}^T \tilde{Z}\bigr)^+ \bigr) = \rank \bigl(
\tilde{Z}^T \tilde{Z} \bigl(\tilde{Z}^T \tilde{Z}\bigr)^+
\bigr) \\
&=& \rank \bigl(\tilde{Z}^T \tilde{Z}\bigr).
\end{eqnarray*}
Furthermore,
\begin{eqnarray*}
\tr \bigl(\tilde{Z} U^T (I-P_\Lambda) U
\tilde{Z}^T \bigr) &=& \tr \bigl(U^T (I-P_\Lambda)
U Z^T (I-P_X) Z \bigr)
\\
& = &\tr \bigl(U^T (I-P_\Lambda) U U^T \Lambda U
\bigr)
\\
&=& \tr \bigl(U^T (I-P_\Lambda) \Lambda U \bigr) = 0,
\end{eqnarray*}
where the last line follows from the fact that $(I-P_\Lambda) \Lambda=
0$. It
follows from (\ref{eqn1}) that
\[
\tr \bigl( (I-P_X) Z Q^{-1} Z^T
(I-P_X) \bigr) \le \operatorname{rank}\bigl(Z^T
(I-P_X) Z\bigr) \sigma_e^{2},
\]
and the second statement has been established. Recall from
Section~\ref{secmain} that $\lambda_{\max}$ is the largest eigenvalue
of $Z^T(I-P_X) Z$, and that $R_i$ is the $q_i \times q$ matrix of 0's
and 1's such that $R_i u = u_i$. Now, fix $i \in\{1,2,\ldots,r\}$
and note that
\[
Q = \bigl(\sigma^2_e\bigr)^{-1}
Z^T (I-P_X) Z + D^{-1} \preceq \bigl(
\sigma^2_e\bigr)^{-1} \lambda_{\max}
I_q + D^{-1}.
\]
It follows that
\[
R_i \bigl( \bigl(\sigma^2_e
\bigr)^{-1} \lambda_{\max} I_q + D^{-1}
\bigr)^{-1} R_i^T \preceq R_i
Q^{-1} R_i^T,
\]
and since these two matrices are both positive definite, we have
\begin{eqnarray*}
\bigl( R_i Q^{-1} R_i^T
\bigr)^{-1} & \preceq & \bigl( R_i \bigl( \bigl(
\sigma^2_e\bigr)^{-1} \lambda_{\max}
I_q + D^{-1} \bigr)^{-1} R_i^T
\bigr)^{-1}
\\
& = & \bigl( \bigl(\sigma^2_e\bigr)^{-1}
\lambda_{\max} + \bigl(\sigma^2_{u_i}
\bigr)^{-1} \bigr) I_{q_i},
\end{eqnarray*}
and this proves that the third statement is true.
\end{pf}

Let $\tilde{z}_i$ and $y_i$ denote the $i$th column of $ \tilde{Z}^T =
((I-P_X) Z)^T$ and the $i$th component of $y$, respectively. Also,
define $K$ to be
\[
\sum_{i=1}^N |y_i| \sqrt
{ \sup_{w \in\mathbb{R}^{N+q}_+} t_i^T \Biggl(
t_i t_i^T + \sum
_{j \in\{1,2,\ldots,N\} \setminus
\{i\}} w_j t_j t_j^T
+ \sum_{j=N+1}^{N+q} w_j
t_jt_j^T+ w_i I_q
\Biggr)^{-2} t_i},
\]
where, for $j=1,2,\ldots,N$, $t_j = \tilde{z}_j$, and for $j \in
\{N+1,\ldots,N+q\}$, the $t_j$ are the standard orthonormal basis
vectors in $\mathbb{R}^q$; that is, $t_{N+l}$ has a one in the $l$th
position and zeros everywhere else.
%
\begin{lemma}
\label{lemmagammafunc}
For any $\sigma^2\in\mathbb{R}_+^{r+1}$, we have
\[
h\bigl(\sigma^2\bigr):= \bigl\| \bigl(\sigma^2_e
\bigr)^{-1} Q^{-1} Z^T (I-P_X) y \bigr\|
\le K<\infty.
\]
\end{lemma}

The following result from \citet{kharhobe2011} will be
used in the proof of Lemma~\ref{lemmagammafunc}.
%
\begin{lemma}
\label{lemKandH}
Fix $n \in\{2,3,\ldots\}$ and $m \in\mathbb{N}$, and let
$x_1,\ldots,x_n$ be vectors in $\mathbb{R}^m$. Then
\[
C_{m,n}(x_1;x_2,\ldots,x_n):=
\sup_{w \in\mathbb{R}^n_+} x_1^T \Biggl( x_1x_1^T
+ \sum_{i=2}^n w_i
x_i x_i^T + w_1 I
\Biggr)^{-2} x_1
\]
is finite.
\end{lemma}
\begin{pf*}{Proof of Lemma~\ref{lemmagammafunc}}
Recall that we defined $\tilde{z}_i$ and $y_i$ to be the $i$th
column of $ \tilde{Z}^T = ((I-P_X) Z)^T$ and the $i$th component of
$y$, respectively. Now,
\begin{eqnarray*}
h\bigl(\sigma^2\bigr) & = & \bigl\| \bigl( Z^T
(I-P_X) Z + \sigma_e^2 D^{-1}
\bigr)^{-1} Z^T (I-P_X) y \bigr\|
\\
& = & \Biggl\| \sum_{i=1}^N \bigl(
\tilde{Z}^T \tilde{Z} + \sigma_e^2
D^{-1} \bigr)^{-1} \tilde{z}_i y_i
\Biggr\|
\\
& \le &\sum_{i=1}^N \bigl\| \bigl(
\tilde{Z}^T \tilde{Z} + \sigma_e^2
D^{-1} \bigr)^{-1} \tilde{z}_i y_i
\bigr\|
\\
& = &\sum_{i=1}^N \Biggl\| \Biggl( \sum
_{j=1}^N \tilde{z}_j
\tilde{z}_j^T + \sigma_e^2
D^{-1} \Biggr)^{-1} \tilde{z}_i y_i
\Biggr\|
\\
& = &\sum_{i=1}^N |y_i|
K_i\bigl(\sigma^2\bigr),
\end{eqnarray*}
where
\[
K_i\bigl(\sigma^2\bigr):= \biggl\| \biggl(
\tilde{z}_i \tilde{z}_i^T + \sum
_{j
\in\{1,2,\ldots,N\} \setminus\{i\}} \tilde{z}_j \tilde{z}_j^T
+ \sigma_e^2 D^{-1} \biggr)^{-1}
\tilde{z}_i \biggr\|.
\]
For each $i \in\{1,2,\ldots, N\}$, define
\[
\hat K_i= \sqrt{\sup_{w \in\mathbb{R}^{N+q}_+} t_i^T
\Biggl( t_i t_i^T + \sum
_{j \in\{1,2,\ldots,N\} \setminus\{i\}} w_j t_j
t_j^T + \sum_{j=N+1}^{N+q}
w_j t_jt_j^T+ w_i
I_q \Biggr)^{-2} t_i},
\]
and notice that $K$ can be written as $K = \sum_{i=1}^N |y_i| \hat
K_i$. Therefore, it is enough to show that, for each $i\in\{1,
2,\ldots, N\}$, $K_i(\sigma^2) \le\hat K_i <\infty$.
Now,
\begin{eqnarray*}
K^2_i\bigl(\sigma^2\bigr) & = &
\tilde{z}_i^T \biggl( \tilde{z}_i
\tilde{z}_i^T + \sum_{j \in\{1,2,\ldots,N\} \setminus\{i\}}
\tilde{z}_j \tilde{z}_j^T +
\sigma_e^2 D^{-1} \biggr)^{-2}
\tilde{z}_i
\\
& = &\tilde{z}_i^T \biggl( \tilde{z}_i
\tilde{z}_i^T + \sum_{j \in
\{1,2,\ldots,N\} \setminus\{i\}}
\tilde{z}_j \tilde{z}_j^T +
\sigma_e^2 \biggl( D^{-1} - \frac{1}{\sigma^2_\bullet}
I_q \biggr) + \frac{\sigma_e^2}{\sigma^2_\bullet} I_q
\biggr)^{-2} \tilde{z}_i
\\
& \le &\sup_{w \in\mathbb{R}^{N+q}_+} t_i^T \Biggl(
t_i t_i^T + \sum
_{j \in\{1,2,\ldots,N\} \setminus\{i\}} w_j t_j t_j^T
+ \sum_{j=N+1}^{N+q} w_j
t_jt_j^T+ w_i I_q
\Biggr)^{-2} t_i
\\[-2pt]
& = &\hat K_i^2.
\end{eqnarray*}
Finally, an application of Lemma~\ref{lemKandH} shows that $\hat
K_i^2$ is finite, and the proof is complete.
\end{pf*}

Let $\chi^2_k(\mu)$ denote the noncentral chi-square distribution
with $k$ degrees of freedom and noncentrality parameter $\mu$.
%
\begin{lemma}
\label{lemchisq}
If $J \sim\chi^2_{k}(\mu)$ and $\gamma\in(0,k/2)$, then
\[
E\bigl[J^{-\gamma}\bigr] \le\frac{2^{-{\gamma}} \Gamma (
{k}/{2}-\gamma
)}{\Gamma ( {k}/{2}  )}.
\]
\end{lemma}
\begin{pf}
Since $\Gamma(x-\gamma)/\Gamma(x)$ is decreasing for $x>\gamma>0$,
we have
\begin{eqnarray*}
E\bigl[J^{-\gamma}\bigr] & = & \sum_{i=0}^\infty
\frac{\mu^i e^{-\mu}}{i!} \int_{\mathbb{R}_+} x^{-\gamma} \biggl[
\frac{1}{\Gamma (
{k}/{2} + i  ) 2^{{k}/{2}+i}} x^{{k}/{2}+i-1} e^{-{x}/{2}} \biggr] \,dx
\\
& = & 2^{-\gamma} \sum_{i=0}^\infty
\frac{\mu^i e^{-\mu}}{i!} \frac{\Gamma ( {k}/{2} +i
-\gamma
)}{\Gamma ( {k}/{2} + i  )}
\\
& \le & 2^{-\gamma} \frac{\Gamma ( {k}/{2} -\gamma )}{\Gamma (
{k}/{2}  )}.
\end{eqnarray*}
\upqed\end{pf}

\subsection{\texorpdfstring{An upper bound on $E  [\llVert{y-W\theta}\rrVert^2|\sigma^2]$}
{An upper bound on E[||y-W theta||2|sigma 2]}}
\label{appub1}

We remind the reader that $\theta= (\beta^T\enskip u^T)^T$, $W = (X\enskip
Z)$, and that $\pi(\theta|\sigma^2,y)$ is a multivariate normal
density with mean $m$ and covariance matrix $V$. Thus,
%
\begin{equation}
\label{eqcondexpone} E \bigl[ \llVert {y-W\theta}\rrVert^2 |
\sigma^2 \bigr] = \tr\bigl(WVW^T\bigr) + \llVert {y-Wm}
\rrVert^2,
\end{equation}
and we have
%
\begin{eqnarray}
\label{eqtraceterm}
&&
\tr\bigl(WVW^T\bigr)
\nonumber\\
&&\qquad = \sigma^2_e \tr(P_X) + \tr \bigl(
P_X Z Q^{-1} Z^T P_X \bigr) - 2
\tr \bigl( Z Q^{-1} Z^T P_X \bigr) + \tr \bigl(
Z Q^{-1} Z^T \bigr)
\nonumber
\\
&&\qquad = p \sigma^2_e - \tr \bigl( Z Q^{-1}
Z^T P_X \bigr) + \tr \bigl( Z Q^{-1}
Z^T \bigr)
\nonumber
\\
&&\qquad = p \sigma^2_e + \tr \bigl( Z Q^{-1}
Z^T (I-P_X) \bigr)
\\
&&\qquad = p \sigma^2_e + \tr \bigl( (I-P_X) Z
Q^{-1} Z^T (I-P_X) \bigr)
\nonumber
\\
&&\qquad \le p \sigma^2_e + \rank\bigl(Z^T
(I-P_X) Z\bigr) \sigma_e^{2}
\nonumber
\\
&&\qquad = (p+t) \sigma^2_e,
\nonumber
\end{eqnarray}
where the inequality is an application of
Lemma~\ref{lemmatracebounds}. Finally, a simple calculation shows
that
\[
y - Wm = (I-P_X) \bigl[ I - \bigl(\sigma^2_e
\bigr)^{-1} Z Q^{-1} Z^T (I-P_X)
\bigr] y.
\]
Hence,
%
\begin{eqnarray}
\label{eqesqterm} \llVert {y - Wm}\rrVert & = & \bigl\llVert {(I-P_X)
y - \bigl(\sigma^2_e\bigr)^{-1}
(I-P_X) Z Q^{-1} Z^T (I-P_X) y}
\bigr\rrVert
\nonumber
\\
& \le &\bigl\llVert {(I-P_X) y}\bigr\rrVert + \bigl\llVert {\bigl(
\sigma^2_e\bigr)^{-1} (I-P_X) Z
Q^{-1} Z^T (I-P_X) y}\bigr\rrVert
\nonumber\\[-8pt]\\[-8pt]
& \le &\bigl\llVert {(I-P_X) y}\bigr\rrVert + \bigl\llVert
{(I-P_X) Z}\bigr\rrVert \bigl\llVert {\bigl(\sigma^2_e
\bigr)^{-1} Q^{-1} Z^T (I-P_X) y}
\bigr\rrVert
\nonumber
\\
& \le &\bigl\llVert {(I-P_X) y}\bigr\rrVert + \bigl\llVert
{(I-P_X) Z}\bigr\rrVert K,\nonumber
\end{eqnarray}
where \mbox{$\llVert \cdot\rrVert $} denotes the Frobenius norm and
the last
inequality uses Lemma~\ref{lemmagammafunc}. Finally, combining
(\ref{eqcondexpone}), (\ref{eqtraceterm}) and (\ref{eqesqterm})
yields
\[
E \bigl[ \llVert {y-W\theta}\rrVert^2 | \sigma^2 \bigr]
\le(p+t) \sigma_e^2 + \bigl(\bigl\llVert
{(I-P_X) y}\bigr\rrVert + \bigl\llVert {(I-P_X) Z}\bigr
\rrVert K \bigr)^2.
\]

\subsection{\texorpdfstring{An upper bound on $E[\llVert{u_i}\rrVert^2|\sigma^2]$}
{An upper bound on E[||u i||2|sigma 2]}}
\label{appub2}

Note that
%
\begin{equation}
\label{eqcondexptwo} E \bigl[ \llVert {u_i}\rrVert^2 |
\sigma^2 \bigr] = E \bigl[ \llVert {R_i u}
\rrVert^2 | \sigma^2 \bigr] = \tr \bigl( R_i
Q^{-1} R_i^T \bigr) + \bigl\| E \bigl[
R_i u | \sigma^2 \bigr] \bigr\|^2.
\end{equation}
By Lemma~\ref{lemmatracebounds}, we have
%
\begin{eqnarray}
\label{eqtracetermtwo} \tr \bigl( R_i Q^{-1}
R_i^T \bigr) & \le & \tr \bigl( R_i
\bigl(Z^T (I-P_X) Z\bigr)^+ R_i^T
\bigr) \sigma_e^2
\nonumber
\\
&&{}+ \tr \bigl( R_i (I-P_{Z^T(I-P_X)Z}) R_i^T
\bigr) \sum_{j=1}^r \sigma^2_{u_j}
\\
& = &\xi_i \sigma_e^2 +
\zeta_i \sum_{j=1}^r
\sigma^2_{u_j}.\nonumber
\end{eqnarray}
Now, by Lemma~\ref{lemmagammafunc},
%
\begin{equation}
\label{eqesqtermtwo} \bigl\| E \bigl[ R_i u | \sigma^2
\bigr] \bigr\| \le\llVert {R_i}\rrVert \bigl\| E \bigl[u | \sigma^2
\bigr] \bigr\| = \llVert {R_i}\rrVert h\bigl(\sigma^2\bigr)
\le\llVert {R_i}\rrVert K.
\end{equation}
Combining (\ref{eqcondexptwo}), (\ref{eqtracetermtwo}) and
(\ref{eqesqtermtwo}) yields
\[
E \bigl[ \llVert {u_i}\rrVert^2 | \sigma^2
\bigr] \le\xi_i \sigma_e^2 +
\zeta_i \sum_{j=1}^r
\sigma^2_{u_j} + \bigl(\llVert {R_i}\rrVert K
\bigr)^2.
\]

\subsection{\texorpdfstring{An upper bound on $E[(\llVert{u_i}\rrVert^2)^{-c}|\sigma^2]$}{An upper bound on E[(||u i||2) -c|sigma 2]}}
\label{appub3}

Fix\vspace*{1pt} $i \in\{1,2,\ldots,r\}$. Given $\sigma^2$, $(R_i Q^{-1}
R^T_i)^{-1/2} u_i$ has a multivariate normal distribution with
identity covariance matrix. It follows that, conditional on
$\sigma^2$, the distribution of $u_i^T (R_i Q^{-1} R^T_i)^{-1} u_i$ is
$\chi^2_{q_i}(w)$. It follows from Lemma~\ref{lemchisq} that, as
long as $c \in(0,1/2)$, we have
\[
E \bigl[ \bigl[u_i^T \bigl(R_i
Q^{-1} R^T_i\bigr)^{-1}
u_i \bigr]^{-c} | \sigma^2 \bigr]
\le2^{-c} \frac{\Gamma ( {q_i}/{2} - c
)}{\Gamma ( {q_i}/{2}  )}.
\]
Now, by Lemma~\ref{lemmatracebounds},
\begin{eqnarray*}
&&
E \bigl[ \bigl( \llVert {u_i}\rrVert^2
\bigr)^{-c} | \sigma^2 \bigr]
\\
&&\qquad = \bigl( \bigl(\sigma^2_e\bigr)^{-1}
\lambda_{\max} + \bigl(\sigma^2_{u_i}
\bigr)^{-1} \bigr)^c E \bigl[ \bigl[u_i^T
\bigl( \bigl(\sigma^2_e\bigr)^{-1}
\lambda_{\max
} + \bigl(\sigma^2_{u_i}
\bigr)^{-1} \bigr) I_{q_i} u_i
\bigr]^{-c} | \sigma^2 \bigr]
\\
&&\qquad \le \bigl( \bigl(\sigma^2_e\bigr)^{-1}
\lambda_{\max} + \bigl(\sigma^2_{u_i}
\bigr)^{-1} \bigr)^c E \bigl[ \bigl[u_i^T
\bigl(R_i Q^{-1} R^T_i
\bigr)^{-1} u_i \bigr]^{-c} |
\sigma^2 \bigr]
\\
&&\qquad \le \bigl( \bigl(\sigma^2_e\bigr)^{-1}
\lambda_{\max} + \bigl(\sigma^2_{u_i}
\bigr)^{-1} \bigr)^c 2^{-c} \frac{\Gamma ( {q_i}/{2} - c  )}{\Gamma (
{q_i}/{2}  )}
\\
&&\qquad \le2^{-c} \frac{\Gamma (
{q_i}/{2} - c  )}{\Gamma ( {q_i}/{2}  )} \bigl[ \lambda^c_{\max}
\bigl(\sigma^2_e \bigr)^{-c} + \bigl(
\sigma^2_{u_i} \bigr)^{-c} \bigr].
\end{eqnarray*}
\end{appendix}

\section*{Acknowledgments}

The authors thank three anonymous reviewers for
helpful comments and suggestions.



\printaddresses

\end{document}